\documentclass[10pt]{amsart}
\usepackage{amssymb}
\usepackage[OT2,OT1]{fontenc}

\let\newpf\proof \let\proof\relax 
\newenvironment{pf}{\newpf[\proofname]}{\qed\endtrivlist}

\def\H{\mathbb {H}}

\def\SL{\mathrm {SL}}
\def\SO{\mathrm {SO}}

\def\be{\begin{equation}}
\def\ee{\end{equation}}

\def\bm{\begin{pmatrix}}
\def\em{\end{pmatrix}}

\def\u{{\mathbb U}}

\def\0{{\mathbf 0}}

\def\cal{\mathcal}

\newcommand{\la}{\label}
\newtheorem{thm}{Theorem}[section]
\newtheorem*{mainthm}{Main Theorem}
\newtheorem{cor}[thm]{Corollary}

\newtheorem{lemma}[thm]{Lemma}

\newtheorem{prop}[thm]{Proposition}

\theoremstyle{remark}
\newtheorem{rem}{Remark}[section]

\numberwithin{equation}{section}

\def\theequation {\thesection.\arabic{equation}}

\def \bn {\hfill \\ \smallskip\noindent}

\theoremstyle{definition}

\newtheorem{definition}{Definition}[section]
\def\proof{\bn {\bf Proof.} }

\def\note#1
{\marginpar
{\tiny $\leftarrow$
\par
\hfuzz=20pt \hbadness=9000 \hyphenpenalty=-100 \exhyphenpenalty=-100
\pretolerance=-1 \tolerance=9999 \doublehyphendemerits=-100000
\finalhyphendemerits=-100000 \baselineskip=6pt
#1}\hfuzz=1pt}

\newcommand{\dist}{\operatorname{dist}}

\newcommand{\inter}{\operatorname{int}}
\renewcommand{\mod}{\operatorname{mod}}

\newcommand{\id}{\operatorname{id}}

\newcommand{\eps}{{\epsilon}}
\newcommand{\De}{{\Delta}}

\newcommand{\C}{{\mathbb C}}

\newcommand{\N}{{\mathbb N}}
\newcommand{\Q}{{\mathbb Q}}
\newcommand{\R}{{\mathbb R}}

\newcommand{\Z}{{\mathbb Z}}

\def\B0{{\bold{0}}}

\catcode`\@=12

\def\Empty{}
\newcommand\oplabel[1]{
  \def\OpArg{#1} \ifx \OpArg\Empty {} \else
  	\label{#1}
  \fi}
		
\newcommand{\comm}[1]{}
\newcommand{\comment}[1]{}

\begin{document}

\bigskip\bigskip

\title{The Ten Martini Problem}
\author{Artur Avila}
\author[Svetlana Jitomirskaya] {Svetlana Jitomirskaya$^\dag$}
\address{
Laboratoire de Probabilit\'es et Mod\`eles al\'eatoires\\
Universit\'e Pierre et Marie Curie--Boite courrier 188\\
75252--Paris Cedex 05, France
}
\email{artur@ccr.jussieu.fr}
\address{
University of California, Irvine, California
}
\email{szhitomi@uci.edu}

\thanks{$^\dag$This work was supported in part by NSF, grant DMS-0300974.}
\date{\today}

\begin{abstract}
We prove the conjecture (known as the ``Ten Martini Problem'' after
Kac and Simon) that the spectrum of the almost Mathieu operator is a
Cantor set for all non-zero values of the coupling and all irrational
frequencies.
\end{abstract}

\setcounter{tocdepth}{1}

\maketitle

\section{Introduction}

The almost Mathieu operator is the Schr\"odinger operator on $\ell^2(\Z)$,
\be
(H_{\lambda,\alpha,\theta} u)_n=u_{n+1}+u_{n-1}+2 \lambda \cos 2 \pi
(\theta+n\alpha) u_n,
\ee
where $\lambda,\alpha,\theta \in \R$ are parameters (called
the {\it coupling}, {\it frequency}, and {\it phase}, respectively), and one
assumes that $\lambda \neq 0$.  The interest in this particular model is
motivated both by its connections to physics and by a remarkable richness of
the related spectral theory.  This has made the latter
a subject of intense research in
the last three decades (see \cite {L3} for a recent historical account and for the
physics background).  Here we are concerned
with the topological structure of the spectrum.

If $\alpha=\frac {p} {q}$ is rational, it is well known that the spectrum
consists of the union of $q$ intervals called {\it bands},
possibly touching at the endpoints.
In the case of irrational $\alpha$, the spectrum $\Sigma_{\lambda,\alpha}$
(which in this case does not depend on $\theta$) has been conjectured for
a long time to be a Cantor set (see a 1964 paper of Azbel \cite {Az}).
To prove this conjecture has been dubbed The Ten Martini Problem by
Barry Simon, after an offer of Mark Kac in 1981, see
Problem 4 in \cite {Sim2}.
For a history of this problem see \cite {L3}.
Earlier partial results include \cite {BS}, \cite {Sin}, \cite {HS},
\cite {CEY}, \cite {L},
and recent advances include \cite {P} and \cite {AK}.
In this paper, we solve the Ten Martini Problem as stated in \cite {Sim2}.

\begin{mainthm}

The spectrum of the almost Mathieu operator is a Cantor set for all
irrational $\alpha$ and for all $\lambda \neq 0$.

\end{mainthm}

It is important to emphasize that the previous results mentioned above
covered a large set of parameters $(\lambda,\alpha)$,
which is both topologically generic (\cite {BS}), and of full Lebesgue
measure (\cite {P}).  As it often happens in the analysis of quasiperiodic
systems, the ``topologically generic'' behavior is quite distinct from the
``full Lebesgue measure'' behavior, and the narrow set of parameters left
behind does indeed lie in the interface of two distinct regimes.
Furthermore, our analysis seems to indicate an interesting characteristic
of the Ten Martini Problem, that the two regimes do not cover nicely the
parameter space and hence there is a non-empty ``critical region'' of
parameters in between (see Remarks \ref {remark2},
\ref {remark1}, \ref {appro} and the comments after
Theorem \ref {open gaps}).

This is to some degree reflected in the structure of the proof.
While the reasoning outside of the critical region can be made
quite effective, in the sense that one essentially identifies
specific gaps in the spectrum\footnote {Related either to gaps of periodic
approximations or to eigenvalues of a dual almost Mathieu operator.}, in
order to be able to cover the critical region we make use of very indirect
arguments.  As an example, we show that absence of Cantor spectrum enables
us to ``analytically continuate'' effective
solutions of a small divisor problem, and it is the
non-effective solutions thus obtained that can be related to gaps in the
spectrum.

This paper builds on a large theory.  Especially important for us are
\cite {CEY}, \cite {J}, \cite {P}, whose methods we improve, but several
other ingredients are needed (such as Kotani Theory \cite {Sim1},
the recent estimates on Lyapunov exponents of \cite {BJ1}).  An important
new ingredient is the use of analytic continuation techniques in the
study of $m$-functions and in extending the reach of the analysis of
Anderson localization.

\subsection{Strategy}

In this problem, arithmetics of $\alpha$ rules the game.  When
$\alpha$ is not very Liouville, it is reasonable to try to deal with the
small divisors.  When $\alpha$ is not very Diophantine, this does not work
and we deal instead with rational approximation arguments.
Let $\frac {p_n} {q_n}$ be the approximants of $\alpha \in \R \setminus \Q$. 
Let
\be
\beta=\beta(\alpha)=\limsup \frac {\ln q_{n+1}} {q_n}.
\ee
The relation between $e^\beta$ and $\lambda$ will play an important role in
our argument, and will decide whether we approach the problem from the
Diophantine side or from the Liouvillian side.
As discussed before, our analysis indicates that there are parameters
that can not be effectively described from either side,
and it is only through the use
of indirect arguments that we can enlarge artificially the Diophantine and
Liouville regimes to cover all parameters.  It should be noted that even
with such tricks, both sides will just about meet in the middle.

Since $\Sigma_{\lambda,\alpha}=\Sigma_{-\lambda,\alpha}$,
it is enough to assume $\lambda>0$.  It is known that the behavior of the
almost Mathieu operator changes drastically at $\lambda=1$
(``metal-insulator'' transition \cite {J}).
Aubry duality shows that
$\Sigma_{\lambda,\alpha}=\lambda \Sigma_{\lambda^{-1},\alpha}$.  So each
$\lambda \neq 1$ admits two lines of attack, and this will be determinant in
what follows.  The case    
$\lambda=1$ was settled in \cite {AK} (after several partial results 
\cite {AvMS}, \cite {HS}, \cite {L}), but it is also recovered in our
approach.

We will work on $\lambda<1$ when approaching from the
Liouville side.
The approach from the Diophantine side is more delicate.  There are actually
two classical small divisor problems that apply to the study of the almost
Mathieu operator, corresponding to Floquet reducibility (for $\lambda<1$)
and Anderson localization (for $\lambda>1$).  An important point
is to attack both problems simultaneously, mixing the best of each
problem (``soft'' analysis in one case, ``hard'' analysis in the other).

A key idea in this paper is that absence of Cantor spectrum implies
improved regularity of $m$-functions in the regime $0<\lambda \leq 1$.
This is proved by analytic continuation techniques.  The improved regularity
of $m$-functions (which is fictitious, since we will prove Cantor spectrum)
will be used both in the Liouville side and in the
Diophantine side.  In the Liouville side, it will
give improved estimates for the continuity of the spectrum with respect to
the frequency.  In the
Diophantine side, it will allow us to use (again) analytic continuation
techniques to solve some small divisor problems
in some situations which are beyond what is expected to be possible.

\begin{rem} \label {remark2}

Since our approach, designed to overcome the difficulties in the interface
of the Diophantine and Liouville regimes, works equally well for other
ranges of parameters, it will
not be necessary in the proof to precisely delimitate a critical region. 
For the reasons discussed in Remarks \ref {remark1}, \ref {appro} and in the
comments after Theorem \ref {open gaps}, the critical region is believed to
contain the parameters such that $\beta>0$ and
$\beta \leq |\ln \lambda| \leq 2 \beta$, the parameters such that
$\beta=|\ln \lambda|$ (respectively, $2 \beta=|\ln \lambda|$)
being seemingly inaccessible (even after artificial extension) by
the Diophantine method (respectively, Liouville method).  It is reasonable
to expect that something should be different in the indicated critical
region.  For instance, it is the natural place to look for possible
counterexamples to the ``Dry Ten Martini'' conjecture (for a precise
formulation see Section \ref{mainthm}).

\end{rem}

\section{Background}

\subsection{Cocycles, Lyapunov exponents, fibered rotation
number} \label {fibered}

A (one-dimensional quasiperiodic $\SL(2,\R)$) {\it cocycle}
is a pair $(\alpha,A) \in \R \times C^0(\R/\Z,\SL(2,\R))$,
understood as a {\it linear skew-product}:
\begin{align}
(\alpha,A):&\R/\Z \times \R^2 \to \R/\Z \times \R^2\\
\nonumber
&(x,w) \mapsto (x+\alpha,A(x) \cdot w).
\end{align}
For $n \geq 1$, we let
\be
A_n(x)=A(x+(n-1)\alpha) \cdots A(x)
\ee
($\alpha$ is implicit in this notation).

Given two cocycles $(\alpha,A)$ and $(\alpha,A')$, a {\it conjugacy}
between them is a continuous $B:\R/\Z \to \SL(2,\R)$ such that
\be
B(x+\alpha)A(x)B(x)^{-1}=A'(x).
\ee

The {\it Lyapunov exponent} is defined by
\be
\lim \frac {1} {n} \int \ln \|A_n(x)\| dx,
\ee
so $L(\alpha,A) \geq 0$.  It is invariant under conjugacy.

Assume now that $A:\R/\Z \to \SL(2,\R)$ is homotopic to the identity.  Then
there exists $\psi:\R/\Z \times \R/\Z \to \R$ and $u:\R/\Z \times \R/\Z \to
\R^+$ such that
\be
A(x) \cdot \bm \cos 2 \pi y \\ \sin 2 \pi y \em=u(x,y)
\bm \cos 2 \pi (y+\psi(x,y)) \\ \sin 2 \pi (y+\psi(x,y)) \em.
\ee
The function $\psi$ is called a {\it lift} of $A$.  Let $\mu$ be any
probability on $\R/\Z \times \R/\Z$ which is invariant by the continuous
map $T:(x,y) \mapsto (x+\alpha,y+\psi(x,y))$, projecting over Lebesgue
measure on the first coordinate (for instance, take $\mu$ as any
accumulation point of $\frac {1} {n} \sum_{k=0}^{n-1} T_*^k \nu$ where
$\nu$ is Lebesgue measure on $\R/\Z \times \R/\Z$).  Then the number
\be
\rho(\alpha,A)=\int \psi d\mu \mod \Z
\ee
does not depend on the choices of $\psi$ and $\mu$, and is called the
{\it fibered rotation number} of
$(\alpha,A)$, see \cite {JM} and \cite {H}.
It is invariant under conjugacies homotopic to the identity.
It immediately follows from the
definitions that the fibered rotation number is a continuous function of
$(\alpha,A)$.

Notice that if $A,A':\R/\Z \to \SL(2,\R)$ and $B:\R/\Z \to \SL(2,\R)$
are continuous such that $A$ is homotopic to the identity and
$B(x+\alpha)A(x)B(x)^{-1}=A'(x)$, then
$\rho(\alpha,A)=\rho(\alpha,A')-k\alpha$, where $k$ is
such that $x \mapsto B(x)$ is homotopic to $x \mapsto R_{kx}$, where
\be
R_\theta=\bm \cos 2 \pi \theta&-\sin 2 \pi \theta \\ \sin 2 \pi \theta&
\cos 2 \pi \theta \em.
\ee

\subsection{Almost Mathieu cocycles, integrated density of states, spectrum}
\label {ids}

Let
\be \la{am}
S_{\lambda,E}=\bm E-2\lambda \cos 2 \pi x & -1 \\ 1 & 0 \em.
\ee
We call $(\alpha,S_{\lambda,E})$, $\lambda,\alpha,E \in \R$, $\lambda \neq
0$ {\it almost Mathieu cocycles}.  A sequence $(u_n)_{n \in \Z}$
is a formal solution of the eigenvalue equation
$H_{\lambda,\alpha,\theta} u=Eu$ if and only if
$S_{\lambda,E}(\theta+n\alpha) \cdot \bm u_n\\u_{n-1} \em=\bm u_{n+1}\\u_n
\em$.

Let
\be
L_{\lambda,\alpha}(E)=L(\alpha,S_{\lambda,E}).
\ee
It is easy to see that
$\rho(\alpha,S_{\lambda,E})$ admits a determination
$\rho_{\lambda,\alpha}(E) \in [0,1/2]$.  We let
\be
N_{\lambda,\alpha}(E)=1-2\rho_{\lambda,\alpha}(E) \in [0,1].
\ee
It follows that $E \mapsto N_{\lambda,\alpha}(E)$ is a continuous
non-decreasing function.  The function $N$ is the usually defined
{\it integrated density of states} of $H_{\lambda,\alpha,\theta}$ if
$\alpha \in \R \setminus \Q$ (for $\alpha \in \Q$, $N$ is the integral of
the density of states over different $\theta$), see \cite {AS} and \cite
{JM}.  Thus defining
\be
\Sigma_{\lambda,\alpha}=\{E \in \R,\, N_{\lambda,\alpha} \text { is not
constant in a neighborhood of } E\},
\ee
we see that (consistently with the introduction)
$\Sigma_{\lambda,\alpha}$ is the spectrum of $H_{\lambda,\alpha,\theta}$ for
$\alpha \in \R \setminus \Q$ (in this case the spectrum does not depend on
$\theta$), while for $\alpha \in \Q$, $\Sigma_{\lambda,\alpha}$ is the union
of the spectra of $H_{\lambda,\alpha,\theta}$, $\theta \in \R$.
One also has
\be
\Sigma_{\lambda,\alpha} \subset [-2-2|\lambda|,2+2|\lambda|].
\ee

Continuity of the fibered rotation number implies that
$N_{\lambda,\alpha}$ depends continuously on $(\lambda,\alpha)$ on
$L^\infty(\R)$.

It turns out that there is a relation between $N$ and $L$, the {\it Thouless
formula}, see \cite {AS}
\be
L(E)=\int \ln |E-E'| dN(E').
\ee
By the Schwarz reflection principle, if $J \subset \R$ is an
open interval where the Lyapunov exponent vanishes, then $E \mapsto
N_{\lambda,\alpha}(E)$ is an increasing analytic function of
$E \in J$\footnote{Since $N+\frac {iL} {\pi}$
is holomorphic in upper half plane
and real on $J.$ This can also be obtained from the Thouless formula.} (and
obviously $J \subset \Sigma_{\lambda,\alpha}$).

We will use several times the following result
\cite {BJ1}.

\begin{thm}[\cite {BJ1}, Corollary 2] \label {L}

Let $\alpha \in \R \setminus \Q$, $\lambda \neq 0$.
If $E \in \Sigma_{\lambda,\alpha}$ then
\be
L_{\lambda,\alpha}(E)=\max \{0, \ln |\lambda|\}.
\ee

\end{thm}

This result will be mostly important for us for what it says about the range
$0<\lambda \leq 1$ (zero Lyapunov exponent on the spectrum).  It will be
also very minorly used in our proof of localization when $\lambda>1$.

\subsection{Kotani theory} \label {kota}

Recall the usual action of $\SL(2,\C)$ on the Riemann sphere $\overline \C$:
$\bm a&b\\c&d \em \cdot z=\frac {a z+b} {c z+d}$.
We can of course define $\SL(2,\C)$ cocycles as pairs $(\alpha,A) \in \R
\times C^0(\R/\Z,\SL(2,\C))$, but it is convenient to view a
$\SL(2,\C)$ cocycle as acting by Moebius transformations:
\begin{align}
(\alpha,A):&\R/\Z \times \overline \C \to \R/\Z \times \overline \C\\
\nonumber
&(x,z) \mapsto (x+\alpha,A(x) \cdot z).
\end{align}
If one lets $E$ become a complex number in the definition of the almost
Mathieu cocycle, we get a $\SL(2,\C)$ cocycle.

Let $\H$ be the upper half plane.  Fix $(\lambda,\alpha)$.
It is well known that there exists a continuous function
$m=m_{\lambda,\alpha}:\H \times \R/\Z
\to \H$ such that $S_{\lambda,E}(x) \cdot m(E,x)=m(E,x+\alpha)$, thus
defining an invariant section for the cocycle $(\alpha,S_{\lambda,E})$:
\be \label {minvariantsection}
(\alpha,S_{\lambda,E}) (x,m(E,x))=(x+\alpha,m(E,x+\alpha)).
\ee
Moreover, $E \mapsto m(E,x)$ is holomorphic on $\H$.

\begin{rem}

In the litterature (for instance, in \cite {Sim2}),
it is more common to find the definition of a pair of
$m$-functions, $m_\pm(x,E)$, which is given in
terms of non-zero solutions $(u_\pm(n))_{n \in \Z}$ of
$H_{\lambda,\alpha,x} u=Eu$
which are $\ell^2$ at $\pm \infty$: $m_\pm(x,E)=-
\frac {u_\pm(\pm 1)} {u_\pm(0)}$.  In this notation we have
$m(x,E)=-\frac {1} {m_-(x,E)}$ (the
relation $S_{\lambda,E}(x) \cdot m(E,x)=m(E,x+\alpha)$ is an immediate
consequence of the definition of $m_-(x,E)$).

\end{rem}

The following result
of Kotani theory \cite {Sim2} will be important in two key parts
of this paper.

\begin{thm} \label {kotani}

Let $\alpha \in \R \setminus \Q$, and assume that $L(\alpha,E)=0$ in an
open interval $J \subset \R$.  Then for
every $x \in \R/\Z$, the functions $E \mapsto m(E,x)$ admit a holomorphic
extension to $\C \setminus (\R \setminus J)$, with values in $\H$.  The
function $m:\C \setminus (\R \setminus J) \times \R/\Z \to \H$ is continuous
in both variables.

\end{thm}

\subsection{Polar sets}

Recall one of the possible definitions of
a polar set in $\C$: it is a set of zero {\it
logarithmic capacity}.  We will need only
some properties of polar sets in $\C$ (see for instance \cite {Ho}):
\begin{enumerate}
\item A countable union of polar sets is polar,
\item The image of a polar set by a non-constant holomorphic
function (defined in some domain of $\C$) is a polar set,
\item Polar sets have Hausdorff dimension zero,
thus their intersections with $\R$ have zero Lebesque measure,
\item Let $U \subset \C$ be a domain and let $f_n:U \to \R$
be a sequence of subharmonic functions which is uniformly bounded in
compacts of $U$.  Then $f:U \to \R$ given by $f=\limsup f_n$ coincides with
its (subharmonic) upper regularization
$f^*:U \to \R$ (given by $f^*(z)=\limsup_{w \to z}
f(w)$) outside a polar set.
\end{enumerate}
We will say that a subset of $\R$ is polar if it is polar as a subset of
$\C$.

The following result on analytic continuation is well known.  We will
quickly go through the proof, since a similar idea will play a role later
in a small divisor problem.

\begin{lemma} \label {anal}

Let $W \subset \C$ be a domain and let $f:W \times \R/\Z \to \C$ be a
continuous function.  If $z \mapsto f(z,w)$ is holomorphic for all $w \in \R/\Z$ and $w
\mapsto f(z,w)$ is analytic for some non-polar set of $z \in W$ then $f$
is analytic.

\end{lemma}

\begin{pf}

We may assume that $|f(z,w)|<1$, $(z,w) \in W \times \R/\Z$.  Let
\be \label {C/n}
f(z,w)=\sum \hat f_z(k) e^{2 \pi i k w}.
\ee
Then $z \mapsto \hat f_z(k)$ is holomorphic and $|\hat f_z(k)|<1$.
Using property (1) of polar sets, 
we obtain that there exists a non-polar set 
$\Delta \subset W$,
$\epsilon>0$, and $k>0$ such that
$|\hat f_z(n)| \leq e^{-\epsilon |n|}$ for $z \in \Delta$ and $|n|>k$.
Let
\be
h(z)=\sup_{|n|>k} \frac {1} {|n|} \ln |\hat f_z(n)|.
\ee
Then, by property (4) of polar sets,
$h^*$ is a non-positive subharmonic function satisfying
$h^*(z) \leq -\epsilon$, $z \in \Delta \setminus X$, where $X$ is polar. 
Since $\Delta$ is non-polar, we conclude that $h^*$ is not identically $0$
in $W$.  It follows from the maximum principle that $h^*(z)<0$, $z \in W$.
Thus for any domain $U \subset W$ compactly contained in $W$,
there exists $\delta=\delta(U)>0$ such that
$h(z) \leq -\delta$, $z \in U$.  We conclude that
\be
\frac {1} {|n|} \ln |\hat f_z(n)| \leq -\delta, \quad
|n|>k, z \in U,
\ee
which implies that (\ref {C/n}) converges uniformly on compacts of
$W \times \{w \in \C/\Z,\, 2 \pi |\Im w|<\delta\}$.
\end{pf}

\section{Regularity of the $m$-functions}

\begin{thm} \label {mregu}

Let $\alpha \in \R \setminus \Q$, $\lambda>0$.  Let $m=m_{\lambda,\alpha}:\H
\times \R/\Z \to \H$ be as in \S \ref {kota}.  Then $m$ is analytic.

\end{thm}

\begin{pf}

Let us show that $m$ has a holomorphic extension to
\be
\Omega_\lambda=\{(E,x),\, \Im E>0,\, 2 \lambda \sinh |2 \pi \Im x|<
\Im E\}.
\ee

We have
\be
S_{\lambda,E}(x) \cdot z=E-2 \lambda \cos (2 \pi x)-\frac {1} {z}.
\ee
For $(E,t)$ satisfying
\be \label {E,t}
\Im E>0, \quad 2 \lambda \sinh |2 \pi t|<\Im E,
\ee
define the half-plane
\be
K^1_{\lambda,E,t}=\{z,\, \Im z>\Im E-2\lambda \sinh |2 \pi t|\} \subset \H,
\ee
and the disk
\be
K^2_{\lambda,E,t}=\left \{|z|<|E|+2 \lambda \cosh |2 \pi t|+\frac
{1} {\Im E-2 \lambda \sinh |2 \pi t|} \right \},
\ee
and let
\be
K_{\lambda,E,t}=K^1_{\lambda,E,t} \cap K^2_{\lambda,E,t},
\ee
which is a domain compactly contained in $\H$ depending continuously on
$(E,t)$ satisfying (\ref {E,t}).
If $(E,x) \in \Omega_\lambda$ then $(E,\Im x)$ satisfies (\ref {E,t}) and
one checks directly that
\be
S_{\lambda,E}(x) \cdot \H \subset K^1_{\lambda,E,\Im x},
\ee
\be
S_{\lambda,E}(x) \cdot K^1_{\lambda,E,\Im x}
\subset K^2_{\lambda,E,\Im x}.
\ee
Since $\Im x=\Im x+\alpha$, we have
\be
S_{\lambda,E}(x+\alpha) \cdot S_{\lambda,E}(x) \cdot
\H \subset K_{\lambda,E,\Im x}.
\ee

Thus, by the Schwarz Lemma applied to $\H$, for every $(E,x) \in \Omega_\lambda$,
\be
S_{\lambda,E}(x-\alpha) \cdots S_{\lambda,E}(x-n\alpha) \cdot
\overline \H
\ee
is a sequence of nested compact sets shrinking to a single point $\hat m(E,x)$.
This implies that $\hat m(E,x)$ is the unique solution to
(\ref{minvariantsection}) in $\H.$
Since $m:\H \times \R/\Z \to \H$ is a continuous function satisfying (\ref
{minvariantsection}),
we must have $\hat m(E,x)=m(E,x)$ for $(E,x) \in \H \times \R/\Z$.

Since holomorphic functions
$m^n:\Omega_\lambda \to \H$ given by
\be
m^n(E,x)=S_{\lambda,E}(x-\alpha) \cdots S_{\lambda,E}(x-n\alpha) \cdot i,
\ee
take values in $\H,$ the sequence $m^n$ is normal. Since it converges pointwise to $\hat m$, we conclude that $\hat m$ is
holomorphic.
\end{pf}

\begin{thm} \label {bjkot}

Let $\alpha \in \R \setminus \Q$, and let $0<\lambda \leq 1$.  Let
$m=m_{\lambda,\alpha}:\H \times \R/\Z \to \H$ be as in \S \ref {kota}.
If $J \subset \Sigma_{\lambda,\alpha}$ is an open interval then $m$ admits
an analytic extension
$m:\C \setminus (\R \setminus J) \times \R/\Z \to \H$.

\end{thm}

\begin{pf}

By Theorems \ref {L} and \ref {kotani}, there exists a continuous
extension $m:\C \setminus (\R \setminus J) \times \R/\Z \to \H$
which is analytic in $E$. 
By Theorem \ref {mregu}, $m$ is also analytic in $x$ for $E \in \H$.
Analiticity in $(E,x)$ then follows by Lemma \ref {anal}.
\end{pf}

\begin{rem}

Notice that the proof of Theorem \ref {mregu} uses strongly that the
dynamics in the basis of the almost Mathieu cocycle
is a rotation (and not, say, a hyperbolic toral automorphism or the skew
shift).  But one may still get weaker results on smoothness of
$m$-functions (in the line of Theorem \ref {bjkot}) for those
dynamics (via estimates in the line of \cite {AK2}).

\end{rem}

\section{Analytic continuation}

\begin{lemma}

Let $\alpha \in \R \setminus \Q$ and let $\phi:\R/\Z \to \R$ be analytic,
and let $\theta=\int_{\R/\Z} \phi(x) dx$.  The following are equivalent:
\begin{enumerate}
\item There exists an analytic function $O:\R/\Z \to \SL(2,\R)$, homotopic
to the identity, such that
\be
O(x+\alpha) R_{\phi(x)} O(x)^{-1}=R_\theta,
\ee
\item There exists an analytic function $\psi:\R/\Z \to \R$ such that
\be \label {psi}
\phi(x)-\theta=\psi(x+\alpha)-\psi(x).
\ee
\end{enumerate}

\end{lemma}

\begin{pf}

Obviously (2) implies (1): it is enough to take $O(x)=R_{-\psi(x)}$.

Let us show that (1) implies (2).
If $O(x) \cdot i=i$ for all $x$ then $O(x) \in \SO(2,\R)$ for all $x$ and
since $O$ is homotopic to the identity we have
$O(x)=R_{-\psi(x)}$ for some analytic function
$\psi:\R/\Z \to \R$ which has to satisfy
$\phi(x)-\theta=\psi(x+\alpha)-\psi(x)$.

Thus we may assume that $O(x_0) \cdot i \neq i$ for some $x_0$.
Notice that
\be
O(x_0+n\alpha) \cdot i=R_{n \theta} O(x_0) \cdot i.
\ee
It follows that if $n_k \alpha \to 0$ in $\R/\Z$ then $2n_k\theta \to 0$ in
$\R/\Z$.  This implies that $\theta=\frac {l} {2} \alpha$ for some $l \in
\Z$.  We have
\be
O(x+\alpha)R_{\phi(x)}O(x)^{-1}=R_{\frac {l} {2} (x+\alpha)} R_{-\frac {l}
{2} x},
\ee
which implies
\be
R_{-\frac {l} {2} (x+\alpha)} O(x+\alpha) R_{\phi(x)}=R_{-\frac {l} {2} x}
O(x),
\ee
and we get
\be
R_{-\frac {l} {2} (x+\alpha)} O(x+\alpha) \cdot i=R_{-\frac {l} {2} x} O(x)
\cdot i.
\ee
 It follows that $R_{-\frac {l} {2} x} O(x) \cdot i=z$ does not depend on
$x$.  Let $Q \in \SL(2,\R)$ be such that $Q \cdot z=i$, and set
\be
S(x)=R_{\frac {l} {2} x} Q R_{-\frac {l} {2} x} O(x).
\ee
Since $O,Q:\R/\Z\to \SL(2,\R),$ where $Q(x)=Q$ are homotopic to the identity, we have that  $S:\R/\Z \to \SL(2,\R)$ is homotopic to the identity and, using that $\theta=\frac{l}2 \alpha,$ we have 
\be
S(x+\alpha)R_{\phi(x)}S(x)^{-1}=R_\theta.
\ee
Moreover, $S(x) \cdot i=i$, so $S(x) \in \SO(2,\R)$ and we have
$S(x)=R_{-\psi(x)}$, $\psi:\R/\Z \to \R$.  It follows that $\psi$ satisfies
(\ref {psi}).
\end{pf}

For $\alpha \in \R \setminus \Q$, $0<\lambda \leq 1$, let
$\Lambda_{\lambda,\alpha}$ be the set of $E$ such that there exists an
analytic function $B_E:\R/\Z \to \SL(2,\R)$, homotopic to the identity, and
$\theta(E) \in \R$, such that
\be
B_E(x+\alpha)S_{\lambda,E}(x)B_E(x)^{-1}=R_{\theta(E)}.
\ee

\begin{thm} \label {beta}

Let $\alpha \in \R \setminus \Q$, $0<\lambda \leq 1$.  Let $J \subset
\Sigma_{\lambda,\alpha}$ be an open interval.  Then
\begin{enumerate}
\item If $\beta=0$ then $\Lambda_{\lambda,\alpha}\cap J=J$,
\item If $\beta<\infty$ then either $\Lambda_{\lambda,\alpha}\cap J$
is polar or $\inter \Lambda_{\lambda,\alpha} \cap J \neq \emptyset$.
\end{enumerate}

\end{thm}

\begin{pf}

Assume that $J \subset \Sigma_{\lambda,\alpha}$ is an open interval.
Let $m=m_{\lambda,\alpha}$ be given
by Theorem \ref {bjkot}, so that
$m:\C \setminus (\R \setminus J) \times
\R/\Z \to \H$ is continuous, $E \mapsto m(E,x)$ is holomorphic and
\be
S_{\lambda,E} \cdot m(E,x)=m(E,x+\alpha).
\ee
Let
\be \label {C_E definition}
C_E(x)=\bm \frac {\Re m(E,x)} {|m(E,x)| (\Im m(E,x))^{1/2}}
& -\frac {|m(E,x)|} {(\Im m(E,x))^{1/2}}\\
\frac {(\Im m(E,x))^{1/2}} {|m(E,x)|}
& 0 \em.
\ee
Then
\be
C_E(x+\alpha) S_{\lambda,E}(x) C_E(x)^{-1} \in \SO(2,\R)
\ee
for $E \in J$, $x \in \R/\Z$.  Since $x \mapsto C_E(x)$ is easily verified to be homotopic
to the identity for $E \in J$, we have
\be
C_E(x+\alpha) S_{\lambda,E}(x) C_E(x)^{-1}=R_{\phi(E,x)}
\ee
for some real-analytic function $\phi:J \times \R/\Z \to \R$.  It follows
that $\phi$ has a holomorphic extension $\phi:Z \to \C$ where
$Z \subset \C \times \C/\Z$ is some domain containing
$J \times \R/\Z$.  So there exists a domain $\Delta \subset \C$
such that $J \subset \Delta$ and
$\Delta \times \R/\Z \subset Z$.  For $E \in \Delta$, let
\be
\phi(E,x)=\sum \hat \phi_E(k) e^{2 \pi i k x}.
\ee

Let $E \in J$ be such that there exists an analytic function
$\psi_E:\R/\Z \to \R$ such that
\be
\int_{\R/\Z} \psi_E(x) dx=0,
\ee
\be
\phi(E,x)-\int_{\R/\Z} \phi(E,x) dx=\psi_E(x+\alpha)-\psi_E(x).
\ee
Then
\be \label {psi1}
\psi_E(x)=\sum \hat \psi_E(k) e^{2 \pi i k x}
\ee
where
\be \label {psi2}
\hat \psi_E(k)=\frac {\hat \phi_E(k)} {e^{2 \pi i k \alpha}-1},
\quad k \neq 0,
\ee
\be \label {psi3}
\hat \psi_E(0)=0.
\ee
We can then define an analytic function $B_E:\R/\Z \to \SL(2,\R)$ by
\be
B_E(x)=R_{-\psi_E(x)} C_E(x),
\ee
which satisfies
\be
B_E(x+\alpha) S_{\lambda,E}(x) B_E(x)^{-1}=R_{\theta(E)}, \quad
\theta(E)=\int_{\R/\Z} \phi(E,x) dx.
\ee

Reciprocally, if there exists an analytic function
$B_E:\R/\Z \to \SL(2,\R)$ homotopic to the identity such that $B_E(x+\alpha)
S_{\lambda,E}(x) B_E(x)^{-1}=R_{\theta(E)}$ for some $\theta(E) \in \R$,
then we can write
\be
O_E(x+\alpha) R_{\phi(x)} O_E(x)^{-1}=R_{\theta(E)},
\ee
where
\be
O_E(x)=B_E(x) C_E(x)^{-1}.
\ee
By the previous lemma, there exists an analytic
function (having average $0$)
$\psi:\R/\Z \to \R$ satisfying $\phi(x)-\int_{\R/\Z} \phi(x)
dx=\psi(x+\alpha)-\psi(x)$.

Notice that
\be
\limsup_{|k| \to \infty} \frac {1} {|k|} \ln \frac {1} {|e^{2 \pi i
k\alpha}-1|}=\beta,
\ee
so that if $\beta=0$ then (\ref {psi1}), (\ref {psi2}), and (\ref {psi3})
really define an analytic function for any $E \in J$, so (1) follows.

Let $a:\Delta \to [-\infty,\beta]$ be given by
\be
a(E)=\limsup_{|k| \to \infty} \frac {1} {|k|} \ln \left |
\frac {\hat \phi_E(k)} {e^{2 \pi i k\alpha}-1} \right |.
\ee
By the previous discussion, $\Lambda_{\lambda,\alpha}=\{E \in J,\,
a(E)<0\}$.  If $\beta<\infty$ then $a$ is $\limsup$ of a
sequence of subharmonic functions which are uniformly
bounded on compacts of $\Delta$.  It follows that $a$ coincides with
its upper regularization
\be
a^*(E)=\limsup_{E' \to E} a(E')
\ee
for $E$ outside some exceptional set which is polar.  Thus the set
$\{E \in J,\, a(E)<0\}$ is either polar (contained in the exceptional set)
or it has non-empty interior.
\end{pf}

\begin{lemma} \label {interior}

Let $\alpha \in \R \setminus \Q$, $\lambda>0$.  Then
$\Lambda_{\lambda,\alpha}$ has empty interior.

\end{lemma}

\begin{pf}

We may assume that $0<\lambda \leq 1$ (otherwise the Lyapunov exponent is
positive on $\Sigma_{\lambda,\alpha}$ which easily implies that
$\Lambda_{\lambda,\alpha}=\emptyset$).
Assume that $J \subset \Lambda_{\lambda,\alpha}$ is an open interval.
Then $J \subset \Sigma_{\lambda,\alpha}$
(since $L_{\lambda,\alpha}(E)=0$ for $E \in J$).  Let $B_E$ be as in the
definition of $\Lambda_{\lambda,\alpha}$.  Then the definition of fibered
rotation number (see \S \ref {fibered}) implies
\be
\rho_{\lambda,\alpha}(E)=\theta(E) (\mod \Z).
\ee
By the analyticity of $\rho$ on $J$ there exists $E \in J$, $l \in \Z$,
such that $\theta(E)=l \alpha (\mod \Z)$.
Let $T_E:\R/\Z \to \SL(2,\R)$ be given by
\be
T_E(x)=R_{-l x} B_E(x).
\ee
Then
\be
T_E(x+\alpha) S_{\lambda,E}(x) T_E(x)^{-1}=\id.
\ee
The conclusion is as in \cite {P}.  For $v \in \R^2$,
\be
S_{\lambda,E}(x) T_E(x)^{-1}v=T_E(x+\alpha)^{-1} v.
\ee
So by (\ref{am}) there exists an analytic $U_v:\R/\Z \to \R$ such that
\be
T_E(x+\alpha)^{-1} \cdot v=\bm U_v(x) \\ U_v(x-\alpha) \em.
\ee
Let
\be
U_v(x)=\sum u^v_n e^{2 \pi inx}.
\ee
It is a standard Aubry duality argument (and can be checked by direct
calculation) that $u_n^v\in\ell^2(\Z)$ is an eigenvector of
$H_{\lambda^{-1},\alpha,0} $ with eigenvalue $\lambda^{-1}E.$  The fact
that we get such an eigenvector for every $v\in\R^2$ contradicts the
simplicity of the point spectrum.
\end{pf}

\begin{rem} \label {blu}

Notice that Lemma \ref {interior} and item (1) of
Theorem \ref {beta} already imply the
Ten Martini Problem in the case $\beta=0$, and we did not need any
localization result (the only recent result we used was Theorem \ref {L}).

\end{rem}

\section{Localization and Cantor spectrum}

We say that the operator $H_{\lambda,\alpha,\theta}$ displays {\it Anderson
localization} if it has pure point spectrum with exponentially decaying
eigenvectors.  This requires $\alpha \in \R \setminus \Q$, and implies
that eigenvalues are dense in $\Sigma_{\lambda,\alpha}$.

\begin{thm} \label {cantor localization}

Let $\alpha \in \R \setminus \Q$, and let $\lambda \geq 1$.  Assume that
$\beta<\infty$.  If $H_{\lambda,\alpha,\theta}$ displays Anderson
localization for a non-polar set
of $\theta \in \R$, then $\Sigma_{\lambda,\alpha}$ is a Cantor set.

\end{thm}

\begin{pf}

Let $\Theta$ be the set of $\theta$ such that $H_{\lambda,\alpha,\theta}$
displays Anderson localization.  If $\theta \in \Theta$, and $E$ is an
eigenvalue for $H_{\lambda,\alpha,\theta}$, let $(u_n)_{n \in \Z}$ be a
non-zero eigenvector.  Then
\be
S_{\lambda^{-1},\lambda^{-1}E} \cdot W(x)=e^{2 \pi i \theta}
W(x+\alpha),
\ee
where
\be
W(x)=\bm U(x) e^{2\pi i \theta}\\
U(x-\alpha)\em,
\ee
and
\be
U(x)=\sum u_n e^{2 \pi i n x}.
\ee
Let $M(x)$ be the matrix with columns $W(x)$ and $\overline {W(x)}$.  Then
\be
S_{\lambda^{-1},\lambda^{-1}E}(x) \cdot M(x)=M(x+\alpha)
\bm e^{2 \pi i \theta}&0\\0&e^{-2 \pi i \theta} \em.
\ee
This implies that $\det M(x)$ is independent of $x$, so
$\det M(x)=c i$ for some $c \in \R$.
Notice that if $c=0$ then
\be
V(x+\alpha)=e^{-4 \pi i \theta} V(x),
\ee
with
\be
V(x)=\frac {U(x)} {\overline {U(x)}}
\ee
(notice that $U(x) \neq 0$ except at finitely many $x$ since $U(x)$ is a
non-constant analytic function)
and in particular, if $n_k \alpha \to 0$ then $2 n_k \theta \to 0$.  So
$2 \theta=k \alpha+l$ for some $k,l \in \Z$.
If $c>0$, we have
\be
S_{\lambda^{-1},\lambda^{-1}E}(x)=Q(x+\alpha) R_\theta Q(x)^{-1}
\ee
where $Q:\R/\Z \to \SL(2,\R)$ is given by
\be
Q(x)=\frac {1} {(2 c)^{1/2}} M(x) \bm 1 & i \\ 1 & -i \em
\ee
and if $c<0$, we have
\be
S_{\lambda^{-1},\lambda^{-1}E}(x)=Q(x+\alpha) R_{-\theta} Q(x)^{-1}
\ee
where $Q:\R/\Z \to \SL(2,\R)$ is given by
\be
Q(x)=\frac {1} {(-2 c)^{1/2}} M(x) \bm 1 & i \\ 1 & -i \em
\bm 1 & 0 \\ 0 & -1 \em.
\ee
It follows that in either case $\lambda^{-1} E
\in \Lambda_{\lambda^{-1},\alpha}$ and moreover,
\be
\rho_{\lambda^{-1},\alpha}(\lambda^{-1}E)=\pm \theta+k \alpha (\mod \Z)
\ee
for some $k \in \Z$.

Let $\Theta' \subset \Theta$ be the set of all $\theta$ such that $2 \theta
\neq k\alpha+l$ for all $k,l \in \Z$.
Let $J \subset \Sigma_{\lambda,\alpha}$ be an open interval.
Then for any $\theta \in \Theta'$, there exists some
$E \in J$ such that $E$ is an eigenvalue for
$H_{\lambda,\alpha,\theta}$, and by the previous discussion any such $E$
satisfies
\be
N_{\lambda^{-1},\alpha}(\lambda^{-1}E)=1-2\rho_{\lambda^{-1},\alpha}
(\lambda^{-1}E)=1-2(\varepsilon \theta+k \alpha+l),
\quad \text {for some } k,l \in \Z, \varepsilon \in \{1,-1\},
\ee
\be
\lambda^{-1}E \in \Lambda_{\lambda^{-1},\alpha}.
\ee
It follows that
\be
\Theta' \subset
\left \{\varepsilon \frac{1-N_{\lambda^{-1},\alpha}
(\Lambda_{\lambda^{-1},\alpha} \cap \lambda^{-1}J)}2-k\alpha-l,\quad
k,l \in \Z, \varepsilon \in \{1,-1\} \right \}.
\ee
By item (2) of
Theorem \ref {beta} and Lemma \ref {interior},
$\Lambda_{\lambda^{-1},\alpha} \cap \lambda^{-1}J$ is polar.
Since $N_{\lambda^{-1},\alpha}$ is a non-constant analytic function on
$\lambda^{-1}J$, it
follows that $\Theta'$ is also polar.  Thus $\Theta \subset \Theta'
\cup \left \{\frac {1} {2} (k\alpha+l),\, k,l \in \Z \right \}$
is polar.
\end{pf}

\begin{rem} \label {remark1}

In \cite {P}, it is shown that if $\alpha \in DC$ then Anderson localization
of $H_{\lambda,\alpha,0}$ implies Cantor spectrum.
We can not however use the argument of Puig
(based on analytic reducibility) 
to conclude Cantor spectrum in the generality we need.
Indeed, we are not able to conclude analytic      
reducibility from localization of $H_{\lambda,\alpha,0}$ in our setting
(in a sense, we spend all our regularity to take care of 
small divisors in the localization result, which is half of analytic
reducibility, and there is nothing left for the other half).  Though this
can be bypassed (using Kotani theory to conclude continuous reducibility
under the assumption of non-Cantor spectrum), there is a much more serious  
difficulty in this approach, see Remark \ref {appro}.

\end{rem}

The next result gives us a large range of $\lambda$ and $\alpha$ where
Theorem~\ref {cantor localization} can be applied.

\begin{thm} \label {localization}

Let $\alpha \in \R \setminus \Q$ be such that $\beta=\beta(\alpha)<\infty$,
and let $\lambda>e^{\frac {16} {9} \beta}$.
Then $H_{\lambda,\alpha,\theta}$ displays Anderson localization for almost
every $\theta$.

\end{thm}

This result improves on \cite {J}, where Anderson localization
was proved under
the assumption that $\alpha$ is Diophantine.  Recall that $\alpha$ is said
to satisfy a Diophantine condition (briefly, $\alpha \in DC$) if
\be
\ln q_{n+1}=O(\ln q_n)
\ee
where $\frac {p_n} {q_n}$ are the rational approximations of $\alpha$.
In particular $\alpha \in DC$ implies (but is strictly stronger than)
$\beta(\alpha)=0$. The proof in \cite{J} with some modifications can be extended to the case $\beta(\alpha)=0$ but not to the case $\beta(\alpha)>0.$

The proof of Theorem~\ref {localization} is the most technical part of this
paper, and the considerations involved are independent from our other
arguments.  We will thus postpone its proof to \S~\ref {proof of
localization}.

\begin{rem} \label {appro}

We expect that the operator $H_{\lambda,\alpha,0}$ does not display Anderson
localization for $1<\lambda \leq e^{2 \beta}$.  The key reason is that in
this regime $0$ is a very resonant phase, and since $\alpha$ is Diophantine 
only in a very weak sense, the compound effect on the small divisors can  
not be compensated by the Lyapunov exponent.  See also Remark~\ref {approx1}.

\end{rem}

\section{Fictitious results on continuity of the spectrum}

The spectrum $\Sigma_{\lambda,\alpha}$ is a continuous function of $\alpha$
in the Hausdorff topology.
There are several results in the literature about quantitative continuity. 
The best general result is due to \cite {AvMS}, $1/2$-H\"older continuity. 
Better estimates can be obtained for $\alpha$ not very Liouville in
the region of positive Lyapunov exponent \cite {JK}.
None of those results are enough for our purposes.

The results described above have something in common: they deal with
something that actually happens, and it is not clear if it is possible to
improve them sufficiently (to the level we need).  Thus we will
argue by contradiction: assuming the spectrum is not Cantor, we will get
very good continuity estimates.  This will allow us to proceed the argument,
but obviously, since we will eventually conclude that the spectrum is a
Cantor set, estimates in this section are not valid for any
existing almost Mathieu operator.  Those estimates might be useful also when
analyzing more general Schr\"odinger operators.

\begin{thm} \label {continuity}

Let $\alpha \in \R \setminus \Q$ and $0<\lambda \leq 1$.  Let
$J \subset \R$ be an open interval such that
$\overline J \subset \inter \Sigma_{\lambda,\alpha}$.  There exists
$K>0$ such that
\be
|N_{\lambda,\alpha}(E)-N_{\lambda,\alpha'}(E)|
\leq K |\alpha-\alpha'|, \quad E \in J.
\ee

\end{thm}

\begin{pf}

Let $m=m_{\lambda,\alpha}$ be as in Theorem \ref {bjkot}.  Define $x \mapsto
C_E(x)$ by (\ref {C_E definition}).
Then, as discussed in the proof of Theorem \ref{beta}, $C_E:\R/\Z \to \SL(2,\R)$ is homotopic to the identity and satisfies
$C_E(x+\alpha)S_{\lambda,E}(x)C_E(x)^{-1} \in \SO(2,\R)$ so
\be
C_E(x+\alpha)S_{\lambda,E}(x)C_E(x)^{-1}=R_{\phi_E(x)},
\ee
where $\phi_E:\R/\Z \to \R$ is analytic.
Recall the definition of the fibered rotation number \S \ref {fibered}.
Then
\be
\rho(\alpha,S_{\lambda,E}(x))=\rho(\alpha,R_{\phi(x)}).
\ee
In this case we can take as lift of $R_{\phi_E(x)}$ the function
$\psi(x,y)=\phi(x)$.

Write
\be
\rho(\alpha',S_{\lambda,E})=\rho(\alpha',C_E(x+\alpha')
S_{\lambda,E}(x)C_E(x)^{-1})=\rho(\alpha',C_E(x+\alpha')C_E(x+\alpha)^{-1}
R_{\phi_E(x)}).
\ee
Since $m$ is analytic in $x,$ we can take as lift of $C_E(x+\alpha')C_E(x+\alpha)^{-1} R_{\phi_E(x)}$
a function $\tilde \psi(x,y)$ satisfying $|\tilde \psi(x,y)-\phi(x)| \leq
K|\alpha-\alpha'|$.  Thus
\be
\|\rho(\alpha,S_{\lambda,E})-\rho(\alpha',S_{\lambda,E})\|_{\R/\Z}
\leq \int \sup_y |\phi(x)-\tilde \psi(x,y)| dx \leq K |\alpha-\alpha'|.
\ee
The result now follows, since $N=1-2\rho$ (see \S \ref {ids})
for the determination of $\rho$ in $[0,1/2]$.
\end{pf}

\begin{rem}

Clearly we also get the fictitious estimate
\be
|L_{\lambda,\alpha'}(E)-L_{\lambda,\alpha}(E)| \leq K|\alpha-\alpha'|, \quad E \in J.
\ee

\end{rem}

\section{Gaps for rational approximants}

It is well known that for any $\lambda \neq 0$ if $\frac {p} {q}$ is a
minimal denomination of a rational number then $\Sigma_{\lambda,\frac {p}
{q}}$ consists of $q$ bands with disjoint interior.  All those bands are
actually disjoint, except if $q$ is even when there are two bands touching
at $0$ \cite {vM}, \cite {CEY}.
The variation of $N_{\lambda,\frac {p} {q}}$
in each band is precisely $1/q$.  The
connected components of
$\R \setminus \Sigma_{\lambda,\frac {p} {q}}$ are called gaps.
Let $M(\lambda,\frac {p} {q})$ be the maximum size of the bands of
$\Sigma_{\lambda,\frac {p} {q}}$.

The following result is well known.

\begin{lemma} \label {density of states}

Let $\alpha \in \R \setminus \Q$ and $\lambda \neq 0$.  If $\frac {p_n}
{q_n} \to \alpha$, then $M(\lambda,p_n/q_n) \to 0$.
In particular (since $N_{\lambda,p_n/q_n} \to N_{\lambda,\alpha}$
uniformly), if one selects a point $a_{n,i}$ in each band of
$\Sigma_{\lambda,p_n/q_n}$ then
\be \la{conv}
\frac {1} {q_n} \sum_i \delta a_{n,i} \to dN_{\lambda,\alpha} \quad \text
{in the weak$^*$ topology.}
\ee

\end{lemma}

In \cite {CEY}, a lower bound for the size of gaps of $\Sigma_{\lambda,\frac
{p} {q}}$
is derived of the form
$C(\lambda)^{-q}$, where, for instance, $C(1)=8$.  We will need the
following sharpening of this estimate, in the case where $\frac {p} {q}$ are
close to a given irrational number.

\begin{thm} \label {cey}

Let $\alpha \in \R \setminus \Q$ and let $0<\lambda \leq 1$.
Let $\frac {p_n} {q_n} \to \alpha$.
For every $\epsilon>0$, for every $n$ sufficiently large,
all gaps of $\Sigma_{\lambda,\frac {p_n} {q_n}}$ have size ar least
$e^{-\epsilon q_n} \lambda^{q_n/2}$.

\end{thm}

\begin{pf}

It is known (see the proof of \cite {CEY} Theorem 3.3 for the
case $\lambda=1$, the general case being obtained as described in the
proof of \cite {CEY} Corollary 3.4) that for any bounded gap
$G$ of $\Sigma_{\lambda,\frac {p} {q}}$, one can find a
sequence $a_i$, $1 \leq i \leq q$, with one $a_i$ in each band of
$\Sigma_{\lambda,\frac {p} {q}}$, such that $G=(a_i,a_{i+1})$ and
\be
\prod_{j \neq i} |a_j-a_i| \geq \lambda^m,
\ee
where $q=2m+1$ or $q=2m+2$.

Let $G_n$ be a bounded gap of $\Sigma_{\lambda,\frac {p_n} {q_n}}$
of minimal size.
Then
\be
|G_n| \geq \lambda^{q_n/2} \prod_{j \neq i_n,i_n+1}
|a_{n,j}-a_{n,i_n}|^{-1},
\ee
where the $a_{n,i}$ satisfy the hypothesis of the previous lemma.
Passing to a subsequence, we may assume that $a_{n,i_n} \to E \in
\Sigma_{\lambda,\alpha}$ and
$|G_n| \to 0$ (otherwise the result is obvious).
By the previous lemma,
we get the estimate for $0<\delta<1$ and for $n$ large
\be
\frac {1} {q_n} \ln (|G_n| \lambda^{-q_n/2}) \geq -\frac {1} {q_n}
\sum_{j \neq i_n,i_n+1} \ln |a_{n,j}-a_{n,i_n}| \geq -
\frac {1} {q_n} \sum_{|a_{n,j}-a_{n,i_n}|>\delta} \ln |a_{n,j}-a_{n,i_n}|,
\ee
which implies by
(\ref{conv}) and the definition of the weak$^*$ topology that
\be
\liminf \frac {1} {q_n} \ln (|G_n| \lambda^{-q_n/2}) \geq
-\int_{|E'-E|>\delta} \ln |E-E'| dN_{\lambda,\alpha}(E').
\ee
Thus
\be
\liminf \frac {1} {q_n} \ln (|G_n| \lambda^{-q_n/2}) \geq
-\int \ln |E-E'| dN_{\lambda,\alpha}(E').
\ee
By the Thouless formula and Theorem \ref {L}, this gives
$\liminf \frac {1} {q_n} \ln (|G_n| \lambda^{-q_n/2}) \geq
-L_{\lambda,\alpha}(E)=0$.
\end{pf}

\begin{rem}

It is possible to get an estimate on the convergence rate on Lemma \ref
{density of states} using
\cite {AvMS}.  This implies an estimate on the rate of convergence in
Theorem \ref {cey}.

\end{rem}

\section{Proof of the Main Theorem}\la{mainthm}

We now put together the results of the previous sections.  Recall that it is
enough to consider $\lambda>0$, and that the case $\lambda=1$ follows from
Theorem 1.5 of \cite {AK}.  Moreover, Cantor spectrum for $\lambda$ implies
Cantor spectrum for $\frac {1} {\lambda}$.  Let $\beta=\beta(\alpha)$.
The Main Theorem follows then from the following.

\begin{thm} \label {1}

Let $\alpha \in \R \setminus \Q$.  Then
\begin{enumerate}
\item If $\beta<\infty$ and
$\lambda>e^{\frac {16} {9} \beta}$,
$\Sigma_{\lambda,\alpha}$ is a Cantor set,
\item If $\beta=\infty$ or if $0<\beta<\infty$ and
$e^{-2\beta}<\lambda \leq 1$, $\Sigma_{\lambda,\alpha}$ is a Cantor set.
\end{enumerate}

\end{thm}

\begin{pf}

Item (1) follows from Theorems \ref {localization} and \ref {cantor
localization}.

To get item (2), we argue by contradiction.  Let $J \subset \inter \Sigma$
be a compact interval.  Then the density of states satisfies $\frac {dN}
{dE} \geq c>0$ for $E \in J$ \cite {AS}.
Let $\frac {p} {q}$ be close to
$\alpha$ such that $\frac {1} {q}
\ln |\alpha-\frac {p} {q}|$ is close to $-\beta$.
By Lemma \ref{density of states} and Theorem \ref {cey},
$J \setminus \Sigma_{p/q}$ contains an interval $G=(a,b)$ of size
$e^{-\epsilon q} \lambda^{q/2}$.  Notice that
$N_{\lambda,\frac {p} {q}}(a)=N_{\lambda,\frac {p} {q}}(b)$.
Theorem \ref {continuity}
implies
\be
|N_{\lambda,\alpha}(a)-N_{\lambda,\alpha}(b)| \leq K \left |\alpha-
\frac {p} {q} \right | \leq e^{\epsilon q} e^{-\beta q}.
\ee
Thus
\be
c \leq \frac {N_{\lambda,\alpha}(a)-N_{\lambda,\alpha}(b)} {a-b} \leq
e^{2 \epsilon q} e^{-\beta q} \lambda^{-q/2}.
\ee
By taking $\epsilon \to 0$, we conclude that $\lambda \leq e^{-2\beta}$.
\end{pf}

Let us point out that $1/2$-H\"older continuity of the spectrum \cite {AvMS}
(which holds for every $\alpha$ and $\lambda$) together with Theorem \ref
{cey} implies the following improvement of \cite {CEY}.  Let us say that
all gaps of $\Sigma_{\lambda,\alpha}$ are open if whenever $E \in
\Sigma_{\lambda,\alpha}$ is such that
$N_{\lambda,\alpha}(E)=k\alpha+l$ for some $k \in \Z \setminus
\{0\}$, $l \in \Z$ then $E$ is the endpoint of some bounded gap (this
obviously implies Cantor spectrum).  The conjecture that
$\Sigma_{\lambda,\alpha}$ has all gaps open for all $\lambda \neq 0$,
$\alpha \in \R \setminus \Q$ is sometimes called the ``dry'' version of the
Ten Martini Problem.

\begin{thm} \label {open gaps}

Let $\alpha \in \R \setminus \Q$ and let $\beta=\beta(\alpha)$.
If $\beta=\infty$ or if $0<\beta<\infty$ and
$e^{-\beta}<\lambda<e^\beta$, $\Sigma_{\lambda,\alpha}$ has all gaps
open.

\end{thm}

The conclusion from Theorem \ref {open gaps}
appears to be the natural boundary of what can be taken honestly from
the Liouvillian method: our computations indicate that although one
can get improved estimates on continuity of the spectrum for
$\lambda>e^\beta$ (following \cite {JK}), things seem to
break up at the precise parameter $\lambda=e^\beta$.  Notice that
$\lambda=e^\beta$ is the expected threshold for localization
(for almost every phase)\footnote{In particular, by the Gordon's argument
enhanced with the Theorem \ref{L}, $H_{\lambda,\alpha,\theta}$ has no 
eigenvalues for $\lambda<e^\beta.$, and no localized eigenfunctions for
$\lambda=e^\beta$.} and 
falls
short of the expected threshold for localization with phase
$\theta=0$, $\lambda=e^{2\beta}$.  Thus the use of fictitious estimates does
not seem to be an artifact of our estimates, but a
rather essential aspect of an
approach that tries to cover all parameters with Diophantine and
Liouvillian techniques.

\begin{rem}

Notice that we do not actually need the measure-theoretical result of \cite
{AK} to obtain Cantor spectrum for $|\lambda|=1$.  Indeed, Lemma \ref{interior} 
and item (1) of Theorem \ref{beta} imply Cantor spectrum for $\beta=0$ (any $\lambda \neq 0$; see
Remark \ref {blu})
and item (2) of Theorem \ref {1} implies
Cantor spectrum for $\beta>0$ (if $|\lambda|=1$).

\end{rem}

\newcommand{\sss}{\setcounter{equation}{0}}

\newcommand{\lra}{\leftrightarrow}
\newcommand{\beq}{\begin{eqnarray}}
\newcommand{\eeq}{\end{eqnarray}}
\newcommand{\bt}{\begin{theorem}}
\newcommand{\et}{\end{theorem}}
\newcommand{\bl}{\begin{lemma}}
\newcommand{\el}{\end{lemma}}
\newcommand{\bc}{\begin{corollary}}
\newcommand{\ec}{\end{corollary}}
\newcommand{\bp}{\begin{prop}}
\newcommand{\ep}{\end{prop}}
\newcommand{\ba}{\begin{array}}
\newcommand{\ea}{\end{array}}
\newcommand{\bs}{\backslash}

\newcommand{\ci}{\cite}
\newtheorem{theorem}{THEOREM}
\renewcommand{\theequation}{\arabic{section}.\arabic{equation}}

\newcommand{\al}{\alpha}
\newcommand{\ga}{\gamma}
\newcommand{\Ga}{\Gamma}
\newcommand{\ch}{\chi}
\newcommand{\ka}{\kappa}
\newcommand{\lb}{\lambda}
\newcommand{\ze}{\zeta}
\newcommand{\Th}{\Theta}
\newcommand{\vp}{\varphi}
\newcommand{\ph}{\phi}
\newcommand{\ps}{\psi}
\newcommand{\Ph}{\Phi}
\newcommand{\ve}{\varepsilon }

\newcommand{\e}{\rm e}
\newcommand{\I}{{\cal I}}
\newcommand{\cH}{\cal H}
\newcommand{\K}{{\cal K}}
\newcommand{\cB}{{\cal B}}
\newcommand{\cV}{{\cal V}}
\newcommand{\cW}{{\cal W}}
\newcommand{\cL}{{\cal L}}

\newcommand{\bi}{\bibitem}

\newfont{\msbm}{msbm10 scaled\magstep1}
\newfont{\msbms}{msbm7 scaled\magstep1} 

\newcommand{\bbq}{\mbox{$\mbox{\msbm Q}$}}
\newcommand{\bbr}{\mbox{$\mbox{\msbm R}$}}
\newcommand{\bbn}{\mbox{$\mbox{\msbm N}$}}
\newcommand{\bbi}{\mbox{$\mbox{\msbm I}$}}
\newcommand{\bbc}{\mbox{$\mbox{\msbm C}$}}
\newcommand{\bbk}{\mbox{$\mbox{\msbm K}$}}

\newcommand{\bbe}{\mbox{$\mbox{\msbm E}$}}
\newcommand{\bbz}{\mbox{$\mbox{\msbm Z}$}}
\newcommand{\bbp}{\mbox{$\mbox{\msbm P}$}}
\newcommand{\bbt}{\mbox{$\mbox{\msbm T}$}}

\newcommand{\bbrs}{\mbox{$\mbox{\msbms R}$}} 
\newcommand{\bbns}{\mbox{$\mbox{\msbms N}$}}
\newcommand{\bbis}{\mbox{$\mbox{\msbms I}$}}
\newcommand{\bbcs}{\mbox{$\mbox{\msbms C}$}}
\newcommand{\bbks}{\mbox{$\mbox{\msbms K}$}}
\newcommand{\bbes}{\mbox{$\mbox{\msbms E}$}}
\newcommand{\bbzs}{\mbox{$\mbox{\msbms Z}$}}
\newcommand{\bbps}{\mbox{$\mbox{\msbms P}$}}

\newcommand{\Bs}{\mbox{$\mbox{\msbms B}$}} 
\newcommand{\Gd}{$G_\delta$}
\newcommand{\Fs}{$F_\sigma$}

\section{Proof of Theorem \ref {localization}} \label {proof of
localization}

We will actually prove a slightly more precise version of Theorem~\ref
{localization}.  Let
\be\la{Theta}\Theta = \{ \theta :\;
 |\sin 2\pi (\theta +(k/2)\al )|< k^{-2}\;
 \mbox{holds for 
infinitely many}\; k'\mbox{s}\}\cup\{\frac{s\pi\al}2,\,s\in\bbz\}.\ee
  $\Th$ is easily seen to have
zero Lebesgue measure by the Borel-Cantelli Lemma.

\begin{thm} \la {loc}

Let $\al\in \bbr\backslash \bbq$ be such that
$\beta=\beta(\al)<\infty,$ and let $\lambda>e^{\frac{16\beta}{9}}.$
Then for $\theta\notin\Th,$ $H_{\lambda,\al,\theta}$ displays Anderson
localization. 

\end{thm}

\begin{rem} \label {approx1}

For $\beta=0$ the theorem holds as well for $\theta=\frac{s\pi\al}2,$
however the proof as presented here will not work.  See \ci{jks} for 
the detail of the argument needed for this case.
In general, we believe that for $\theta$ of the form
$\frac{s\pi\al}2,\,s\in\bbz$ the localization would only hold for
$\lambda>e^{2\beta}$.

\end{rem}

\begin{rem}

We believe that for $\theta\notin\Th$, localization should hold
for $\lambda>e^{\beta}$.  The proof of this fact would require some
additional arguments.  Moreover, for $\lambda \leq e^\beta$, we do not
expect any exponentially decaying
eigenvectors. 

\end{rem}

\begin{rem}

The bound $k^{-2}$ in (\ref{Theta}) can be replaced by any other
sub-exponential function without significant changes in the proof.

\end{rem}


We will use the general setup of \ci{J}, however our
 key technical procedure will have to be quite different.

A formal solution $\Psi _E(x) $ of the equation
 $H_{\lambda,\al,\theta}\Psi _E = E\Psi _E$ will
be called {\it a generalized eigenfunction} if 
\be \la{ge} |\Psi_E(x)| \le C(1+|x|)\ee
for some $C=C(\Psi _E)<\infty $. The corresponding $E$ is called a
 {\it generalized eigenvalue.}
It is well known that to prove Theorem \ref{loc} it suffices to prove
 that generalized eigenfunctions  decay exponentially \ci{ber}.

We will use the notation $G_{[x_1,x_2]}(x,y)$ for matrix elements of the Green's function
$(H-E)^{-1}$ of the operator
 $H_{\lb,\al,\theta}$ restricted to the interval $[x_1,x_2]$ with zero boundary
 conditions at $x_1-1$ and $x_2+1$.
We now fix $\lb,\al$ as in Theorem \ref{loc}.

Fix a generalized eigenvalue $E$, and let $\Psi$ be the corresponding
generalized eigenfunction.  Then
\be \la{bj}
L(E)=\ln\lambda>0.\ee $\lambda$ will
enter into our analysis through $L$ only and it will be convenient to
use $L$ instead. To simplify 
notations,
 in some cases
the $E,\lb,\al$-dependence of various quantities will be omitted.

%
%
Fix $m>0$.  A point $y \in {\bbz}$
 will be called 
{\it $(m,k)$-regular} if there exists an interval $[x_1,x_2],\;x_2=x_1+k-1,$
containing $y,$ such that
$$
|G_{[x_1,x_2]}(y,x_i)|<e^{- m|y-x_i|},\;\mbox{and}\;\;{\rm dist} (y,x_i) 
\ge \frac {1} {40} k;\; i=1,2.
$$
Otherwise, $y$ will be called {\it $(m,k)$-singular.}



It is well known and can be checked easily that values of any formal 
solution $\Psi$ of the equation
$H\Psi=E\Psi$ at a point $x\in I = [x_1,x_2]\subset \bbz$  can be reconstructed from the 
boundary values via

\be \la{poi}
\Psi (x)=-G_{I}(x,x_1)\Psi (x_1-1)-
 G_{I}(x,x_2)\Psi (x_2+1).
\ee

This implies that if $\Psi_E$ is a generalized eigenfunction, then every 
point $y \in \bbz$
with $\Psi_E(y)\not= 0$ is $(m,k)$-singular for $k$ sufficiently large:
$k> k_1(E,m,\theta ,y).$  We assume without loss of generality that
$\Psi(0)\not=0$ and normalize $\Psi$ so that $\Psi(0)=1.$  Our strategy will be to show first that every
sufficiently large $y$ is $(m,\ell(y))$-regular for appropriate
$(m,\ell).$ While $\ell$ will vary with $y,$ $m$ will have a uniform
lower bound. This will be shown in subsections \ref{ress} and
\ref{nonres}.  Exponential decay will be derived out of this property via a
``patching argument'' in subsection \ref{pat}.

Let us denote
$$
P_k(\theta)=\det \left[ (H_{\lb,\al,\theta}-E)\bigg|_{[0,k-1]}\right].
$$
Then the $k-$step transfer-matrix $A_n(\theta)$ (which is the $k$-th iterate
of the Almost Mathieu cocycle, $A_k(\theta)=S_{\lambda,E}(\theta+(k-1)
\alpha) \cdots S_{\lambda,E}(\theta)$) can be written as
\be \la{matr}
A_k(\theta)=\bm P_k(\theta) &-P_{k-1}(\theta+\al)\\
P_{k-1}(\theta)
&-P_{k-2}(\theta+\al)\em\,.  
\quad \ee

Herman's subharmonicity trick \ci{H} yields $\int_0^1\ln
|P_k(\theta)|d\theta
\ge k\ln\lambda;$
together with (\ref{bj}) this implies that
there exists 
$\theta \in [0,1] \;\mbox{with}\; |P_k(\theta) |\ge e^{kL(E)}$. Note that this is the only place in the proof of localization where we have used (\ref{bj}). While this is not really necessary (the rest of the proof can proceed, with only minor technical changes, under the assumption of the lower bound on only one of the four matrix elements, which follows immediately from the positivity of $L(E)$) it simplifies certain arguments in what follows.

By an application of the Cramer's rule we have that 
for any $x_1,\; x_2=x_1+k-1,\; x_1 \le y \le x_2,$ 
\be \label {cra}
\begin{array}{cc}
& |G_{[x_1,x_2]}(x_1,y)|=\left| \frac { P_{x_2-y}(\theta + (y+1)\al )}
{P_k(\theta +x_1\al )}\right|,\\[0.15 in]
& |G_{[x_1,x_2]}(y,x_2)|=\left| \frac { P_{y-x_1}(\theta  +x_1\al )}
{P_k(\theta +x_1\al )}\right|.
\end{array}
\ee

The numerators in (\ref{cra}) can be bounded uniformly in $\theta$
\ci{J,fur}. Namely, for every 
 $E \in \bbr, \, \eps>0,$
there exists $k_2(\eps, E,\al)$ such that 

\be \la{pol0}
|P_n(\theta)|<e^{(L(E)+\eps)n}\ee
for
all $n> k_2(\eps, E,\al),$ all $\theta.$

$P_k(\theta)$ is an even function of $\theta + \frac {k-1}{2}\al$ and
 can be written as a polynomial of degree $k$ in 
$\cos 2\pi(\theta+\frac {k-1}
{2}\al):$ $$ P_k(\theta)
=\sum_{j=0}^k
c_j\cos ^j 2\pi(\theta + \frac {k-1}{2}\al)\stackrel{\mathrm{def}}{=}Q_k(\cos 
2\pi(\theta + \frac {k-1}{2}\alpha)).$$


Let $A_{k,r}=\{\theta \in \R,\, |Q_k(\cos 2\pi\theta)|\le e^{(k+1)r}\}$.
The next lemma shows that every singular point ``produces" a long piece of the
trajectory of the rotation consisting of points 
belonging to an appropriate
 $A_{k,r}$.

\bl \la {prop} Suppose $y\;\in \; {\bbz}$ is
$(L-\eps,k)$-singular, $0<\epsilon<L$.
Then for any $\eps_1>0,\; \frac {1} {40} \leq \delta<1/2,$ for
sufficiently large $k>k(\eps,E,\al,\eps_1,\delta),$
 and for any $x \in \bbz$ such that
$y-(1-\delta) k\le x \le y-\delta k,$
we have that $\theta+(x+ \frac{k-1}{2})\al$  belongs to $A_{k,L-\eps\delta+\eps_1}.$
\el

\begin{pf}

Follows immediately from the definition of regularity,
(\ref{cra}), and (\ref{pol0}).
\end{pf}

\comm{
Lemma \ref {prop} will be always applied to exploit singularity of zero as
follows.

\begin{cor}

Let $\theta_j=\theta+j\alpha$.  For any $\epsilon>0$, for sufficiently large
$k=k(E,\alpha,\epsilon)$, we have
\begin{enumerate}
\item For any $x \in [-k+[\frac {k} {20}],k-[\frac {k} {20}]]$,
$\theta_j \in A_{2k-1,L-\epsilon}$.
\item For any $\frac {1} {2}<\rho \leq \frac {39} {40}$, for any
$x \in [(1-2 \rho) k,(2 \rho-1)k]$, $\theta_j \in A_{2k-1,\rho L+\epsilon}$.

\end{cor}
}

The idea now is to show that $A_{k,r}$ cannot contain $k+1$ uniformly
distributed points. In order to quantify this concept of uniformity we
introduce the following.

\begin{definition}

We will say that the set $\{\theta_1,\ldots,\theta_{k+1}\}$ is
{\it $\eps$-uniform} if 

\be \la{prod}
\max_{z\in[-1,1]}\max_{j=1,\ldots,k+1}\prod_{^{\ell=1}_{\ell \not= j}}^{k+1}\frac{
|z-\cos 2\pi\theta_\ell)|}
{|\cos 2\pi\theta_j-\cos 2\pi\theta_\ell)|}
< e^{{k\eps}}
\ee

\end{definition}

Note that we will use this terminology with ``large'' values of $\eps$ as well. $\eps$-uniformity (the smaller $\eps$ the better) involves uniformity along with certain cumulative
repulsion
of $\pm \theta_i$(mod $1$)'s.

\bl \la{nongen}
Let $\eps_1 <\eps.$ If $\theta_1,\ldots,\theta_{k+1}\in A_{k,L-\eps}$ and
$k>k(\eps,\eps_1)$ is sufficiently large, then
$\{\theta_1,\ldots,\theta_{k+1}\}$ is not $\eps_1$-uniform.
\el

\begin{pf}

Write polynomial $Q_k(z)$
 in the Lagrange interpolation form using $\cos 2\pi\theta_1,\ldots,\cos
2\pi\theta_{k+1}:$
\be \la{lagr}
|Q_k(z)|=\left
|\sum_{j=1}^{k+1}Q_k(\cos 2\pi\theta_j)\frac {\prod_{\ell \not= j}
(z-\cos 2\pi\theta_\ell)}
{\prod_{\ell \not= j}(\cos2\pi\theta_j-\cos 2\pi\theta_\ell)} \right |
\ee
 Let $\theta_0$ be such
 that $|P_k(\theta_0)|\ge
 \exp (kL )$.  The lemma now follows immediately from
 (\ref{lagr}) with $z=\cos 2\pi(\theta_0+\frac {k-1}{2}\al)$.
\end{pf}

Suppose we can find two intervals, $I_1$ around $0$ and $I_2$ around $y$,
of combined length $|I_1|+|I_2|=k+1,$\footnote {Here and in what follows,
the ``length'' $|I|$ of an interval $I=[a,b] \subset \Z$ denotes
cardinality, $|I|=b-a+1$.}
such that we can establish the uniformity of
$\{\theta_i\}$ where
$\theta_i=\theta+(x+ \frac{k-1}{2})\al,\; i=1,\ldots,k+1,$ for $x$ ranging through $I_1\cup
I_2.$ Then we can apply 
Lemma \ref{prop} and Lemma \ref{nongen} to show regularity of $y$. 
This is roughly
going to be the framework for our strategy to establish
regularity. The implementation will depend highly on the position of
$k$ with respect to the sequence of denominators $q_n.$ 

Assume without loss of generality that $k>0.$ Define
$b_n=\max\{q_n^{8/9},\frac 1{20}q_{n-1}\}$ Find $n$ such that
 $b_n<k\le b_{n+1}.$
We will distinguish between the two cases:
\begin{enumerate}
\item {\bf Resonant:} meaning $|k-\ell q_n| \le b_n$ for some
  $\ell\ge 1$
   and
\item {\bf Non-resonant:} meaning $|k-\ell q_n| > b_n$ for all
  $\ell\ge 0.$
\end{enumerate}

We will prove the following estimates.

\bl \la{non-res} Assume $\theta\notin\Th$.  Suppose $k$ is non-resonant. 
Let $s\in \bbn\cup \{0\}$ be the
largest number such that 
$sq_{n-1}\le \dist(k, \{\ell q_n\}_{\ell\ge 0}) \equiv k_0.$
Then for any $\eps>0$ for sufficiently large $n$,
\begin{enumerate}
\item If $s\ge 1$ and $L>\beta$, $k$ is $(L-\frac{\ln
  q_n}{q_{n-1}}-\epsilon,2sq_{n-1}-1)$-regular.
\item If $s=0$
then $k$ is either $(L-\eps, 2[\frac {q_{n-1}} {2}]-1)$ or
$(L-\epsilon,2[\frac {q_n} {2}]-1)$ or
$(L-\epsilon,2 q_{n-1}-1)$-regular.
\end{enumerate}
\el

\bl \la{res} Let in addition $L > \frac{16}9\beta.$
Then for sufficiently large $n,$ every resonant $k$ is $(\frac{L}{50},2q_n-1)$-regular.
\el


We will prove Lemma \ref{non-res} in Subsection \ref{nonres} and Lemma
\ref{res} in Subsection \ref{ress}.  Note that these two subsections
are not independent: the proof of Lemma \ref{res} uses a corollary
of the proof of Lemma \ref{non-res} as an important ingredient.
As our proofs rely on establishing
$\eps$-uniformity of certain quasiperiodic sequences, we will use
repeatedly estimates on trigonometric products that we prove in
Subsection \ref{trig}.

Theorem \ref {loc} can be immediately derived from Lemmas \ref{non-res},
\ref{res} via a ``patching argument'' which we describe now.  (A patching
argument will also be used in one step of the proof of Lemma \ref {res}.)

\subsection{Patching.
Proof of Theorem \ref {loc} assuming Lemmas \ref {non-res} and
\ref {res}} \label {pat}
It is an important technical ansatz of the multiscale analysis that the
exponential
decay of a Green's function at a scale $k$ under certain conditions
generates exponential decay with the same rate at a 
larger scale.
The proof is usually done using block-resolvent expansion, with the
combinatorial factor being killed by the growth of scales.
The proof of Theorem \ref{loc} will consist, roughly, of adapting this type
of argument to our situation.

Fix a generalized eigenvalue $E$ of $H_{\lb,\al,\theta}$, and let $\Psi$ be the corresponding
generalized eigenfunction.

Assume without 
loss of generality that $k$ is positive. Find $n$ so that
$k>q_n$.
We assume that $n$ is sufficiently
large. Let $L_1=\frac {L} {50} \leq L-\beta$.
By Lemmas \ref{non-res} and \ref{res} and definition of
regularity, for any $y>b_n$
there exists an interval $y\in I(y)=[x_1,x_2]\subset \bbz$ such that
\be \la{I'}
\dist(y,\partial I(y))>\frac 1{40}|I(y)|,
\ee
\be \label {I''}
|I(y)|\ge q_n^{8/9}-2,
\ee
\be \label {I}
G_{I(y)}(k,x_i) < e^{-L_1|k-x_i|},\;i=1,2.
\ee
In addition, if $b_j<y\le b_{j+1}$ we have 
\be\la{b}|I(y)|\le 2q_j.\ee
We denote the boundary of the interval $I(y)$, the set $\{x_1,x_2\},$ by
 $\partial I(y).$ For $z\in
\partial I(y)$ we let $z^\prime$ be the neighbor of $z,$ (i.e.,
 $|z-z^\prime|=1)$ not belonging to
 $ I(y).$ 

We now expand  $\Psi(x_2+1)$ in (\ref{poi}) iterating  (\ref{poi}) 
with $I= I(x_2+1).$ In case $q_n^{8/9}<x_1-1$ we also expand $\Psi(x_1-1)$
 using (\ref{poi}) with $I=I(x_1-1).$ 
 We continue to expand each term of the
form $\Psi(z)$  in the same fashion until we arrive to $z$ such that either 
$z \leq b_n,\; z>k^2$ or the number of $G_I$ terms in the product
becomes $[\frac{40k}{q_n^{8/9}}],$ whichever comes first. We then obtain an expression 
of the form
\be\la{block}
\Psi(k)=\displaystyle\sum_{s ; z_{i+1}\in\partial I(z_i^\prime)} 
G_{I(k)}(k,z_1) G_{I(z_1^\prime)}
(z_1^\prime,z_2)\cdots G_{I(z_s^\prime)}
(z_s^\prime,z_{s+1})\Psi(z_{s+1}^\prime).
\ee
where in each term of the summation we have $z_i>b_n$, $i=1,\ldots,s,$ and 
either $0<z_{s+1}^\prime \leq b_n\;$,
 $s\le \frac{40k}{q_n^{8/9}}$
 or $z_{s+1}^\prime >k^2\;$,
 $s\le \frac{40k}{q_n^{8/9}}$ or $s+1=[40\frac{k}{q_n^{8/9}}] .$ 
By construction, for  each $z_i^\prime,\, i\le s,$ we have that
$I(z_i^\prime)$ is well-defined and satisfies (\ref{I}), (\ref{b}).  We now 
consider the three cases,
$0<z_{s+1}^\prime \leq b_n,$ $z_{s+1}^\prime >k^2\;$
and $s+1=[\frac{40k}{q_n^{8/9}}]$ separately. 
If $0<z_{s+1}^\prime \leq b_n$ we have, by (\ref{I}) and (\ref{ge}),
\begin{align}
| G_{I(k)}(k,z_1) G_{I(z_1^\prime)}
&(z_1^\prime,z_2)\cdots G_{I(z_s^\prime)}
(z_s^\prime,z_{s+1})\Psi(z_{s+1}^\prime)|\\
\nonumber
&\le C e^{-L_1(|k-z_1|+\sum_{i=1}^{s}|z_i^\prime-z_{i+1}|)}(1+b_n)\\
\nonumber
& \leq C e^{-L_1(|k-z_{s+1}|-(s+1))}(1+b_n)
\leq C e^{-L_1(k-b_n-\frac{40k}{q_n^{8/9}})}(1+b_n).
\end{align}

Similarly, if $z_{s+1}^\prime >k^2,\;$ we use (\ref {I}) and
(\ref{b}) to get
$$
| G_{I(k)}(k,z_1) G_{I(z_1^\prime)}
(z_1^\prime,z_2)\cdots G_{I(z_s^\prime)}
(z_s^\prime,z_{s+1})\Psi(z_{s+1}^\prime)|\le C e^{-L_1(k^2-k
-\frac{40k}{q_n^{8/9}})}(1+3k^{9/4}).
$$   
Finally, if $s+1=[\frac{40k}{q_n^{8/9}}], $ using again 
(\ref{ge}), (\ref{I}), and also (\ref {I'})
 we can estimate
$$| G_{I(k)}(k,z_1) G_{I(z_1^\prime)}
(z_1^\prime,z_2)\cdots G_{I(z_s^\prime)}
(z_s^\prime,z_{s+1})\Psi(z_{s+1}^\prime)|\le C e^{-L_1
\frac {1} {40}q_n^{8/9}[\frac{40k}{q_n^{8/9}}]}(1+k^2) .$$     
 In either case,  
\be\la{est}
 | G_{I(k)}(k,z_1) G_{I(z_1^\prime)}
(z_1^\prime,z_2)\cdots G_{I(z_s^\prime)}
(z_s^\prime,z_{s+1})\Psi(z_{s+1}^\prime)|\le e^{-\frac{9L_1}{10}k}
\ee
for $k$ sufficiently large. Finally, we observe that
the total
 number of terms in
 (\ref{block})
is bounded above by $2^{[\frac{40k}{q_n^{8/9}}]}.$ Combining it with
(\ref{block}),
(\ref{est}) we obtain
$$
 |\Psi(k)|\le 2^{[\frac{40k}{q_n^{8/9}}]}e^{-\frac{9L_1}{10}k}<e^{-\frac{4L_1}{5}k}
$$      
for large $k$.\qed   

  \subsection {Estimates on trigonometric products.}\la{trig}

We will use the notation $\|z\|_{\R/\Z}$ for
the distance to the nearest integer.

\bl 

Let $p,q$ be relatively prime.  We have

\begin{enumerate}

\item 
Let $1\le k_0\le q$ be such that $|\sin
  2\pi(x+\frac{k_0p}{2q})|=\min_{1\le k\le q}|\sin
  2\pi(x+\frac{kp}{2q})|$. Then
\be \la{rat-1}
\ln q +\ln \frac{2}{\pi} <\sum_{^{k=1}_{k\not=k_0}}^q\ln|\sin
  2\pi(x+\frac{kp}{2q})|+(q-1)\ln 2 \le \ln q,
\ee
\item
\be \la{rat0}
\sum_{k=1}^{q-1}\ln |\sin \frac {\pi kp}{q}|=-(q-1)\ln 2 +\ln q.
\ee
\end{enumerate}
\el

\begin{pf}

We use that
\be\la{lnsin}
\ln |\sin \frac x2|=-\ln 2-\frac 12\sum_{k\not=0}\frac{e^{ikx}} {|k|}.
\ee
Thus, for $x\not= \frac{k\pi}{2q},$
\begin{align}
\sum_{j=1}^q\ln|\sin
  2\pi(x+\frac{jp}{2q})|&=-q\ln 2 -\frac 12\sum_{k\not=0}\frac
  {1} {|k|}\sum_{j=1}^qe^{2\pi ik(2x+\frac{jp}q)}\\
\nonumber
&=-q\ln 2-\frac 12\sum_{k\not=0}\frac
  {1} {|k|} e^{4\pi ikqx}=-q\ln 2 +\ln 2+\ln |\sin 2\pi qx|.
\end{align}
Thus 
\be \la{fr}
 \displaystyle\sum_{^{k=1}_{k\not=k_0}}^q\ln|\sin
  2\pi(x+\frac{kp}{2q})|+(q-1)\ln 2 =\ln \frac{|\sin2\pi q(x+\frac{k_0p}{2q})|}{|\sin 2\pi (x+\frac{k_0p}{2q})|}
  \ee
 
It is easily checked that if $0<qx\le\frac{\pi}2$, then
$\frac{2q}{\pi}<\frac{\sin qx}{\sin x}<q$.
Since $\|2 x+\frac{k_0p}{q}\|_{\R/\Z}\le\frac1{2q}$, (\ref{fr}) implies
(\ref{rat-1}).  (\ref{rat0}) follows by taking the limit in (\ref{fr}).
\end{pf}

For $\alpha\notin \bbq$ let $\frac{p_n}{q_n}$ be its continued
fraction approximants. Set $\Delta_n=|q_n\alpha-p_n|$.  We recall the
following basic estimates
\be \la{cf1}
\frac {1} {q_n}>\Delta_{n-1}>\frac {1} {q_n+q_{n-1}},
\ee
\be \la{cf2}
\|k \alpha\|_{\R/\Z} > \Delta_{n-1}, \quad q_{n-1}+1 \leq k \leq q_n-1.
\ee

Notice that if $z,w \in \R$ are such that $\cos(z-w) \geq 0$ then
\be \label {sin z-w}
\left |\frac {\sin z} {\sin w}-1 \right | \leq
\left |\cos (z-w)-1+\frac {\cos w} {\sin w} \sin (z-w) \right | \leq \left
|2 \frac {\sin (z-w)} {\sin w} \right |.
\ee

\bl

%
Let $1\le k_0\le q_n$ be such that $|\sin
  2\pi(x+\frac{k_0\al}{2})|=\displaystyle\min_{1\le k\le q_n}|\sin
  2\pi(x+\frac{k\al}{2})|$. Then
\be \la{irrat-1}
\left |\sum_{^{k=1}_{k\not=k_0}}^{q_n}\ln|\sin
  2\pi(x+\frac{k\al}{2})|+(q_n-1)\ln 2 \right | < C\ln q_n.
\ee
%
\el

\begin{pf}

Let $1\le k_1 \le q_n$ be such that
\be
|\sin 2\pi(x+\frac{k_1 p_n}{2 q_n})|=\min_{1\le k\le q_n}|\sin
  2\pi(x+\frac{k p_n}{2 q_n})|.
\ee
We first remark that, by (\ref{cf2})
\be \label {k-k'}
\|(2x+k\alpha)-(2x+k'\alpha)\|_{\R/\Z}=\|(k-k')\alpha\|_{\R/\Z} \geq
\Delta_{n-1}, \quad 1 \leq k,k' \leq q_n, k \neq k'.
\ee
Applying this to the case $k'=k_0$, we get, by (\ref{cf1}),
$|\ln|\sin 2\pi(x+\frac{k \alpha}{2})||<C \ln q_n$, $k \neq k_0$.  An even
simpler argument
\be \label {k-k''}
\|(2x+k\frac {p_n} {q_n})-(2x+k'\frac {p_n} {q_n})\|_{\R/\Z}=
\|(k-k') \frac {p_n} {q_n}\|_{\R/\Z} \geq \frac {1} {q_n},
\quad 1 \leq k,k' \leq q_n, k \neq k',
\ee
also gives that if $k \neq k_1$ then
$|\ln|\sin 2\pi(x+\frac{k p_n}{2 q_n})||<C \ln q_n$.  This and (\ref {rat-1})
show that it is enough to get the estimate
\be \label {C ln q_n}
\sum_{^{k=1}_{k \neq k_0,k_1}}^{q_n} \ln|\frac{\sin
  2\pi(x+\frac{k\al}{2})}{\sin
  2\pi(x+\frac{kp_n}{2q_n})}|<C \ln q_n.
\ee
By (\ref {sin z-w}),
\be \label {sinka}
|\frac{\sin
  2\pi(x+\frac{k\al}{2})}{\sin
  2\pi(x+\frac{kp_n}{2q_n})}-1|<\frac {C_0 \Delta_n} {|\sin
  2\pi(x+\frac{kp_n}{2q_n})|},
\ee
so we have
\be
\ln|\frac{\sin
  2\pi(x+\frac{k\al}{2})}{\sin
  2\pi(x+\frac{kp_n}{2q_n})}|<\frac {C \Delta_n} {|\sin
  2\pi(x+\frac{kp_n}{2q_n})|},
\ee
provided that $C_0 \Delta_n<\frac 14 |\sin 2\pi(x+\frac{kp_n}{2q_n})|$.
Let $s_1,...,s_r$ be an enumeration of
$\{1 \leq k \leq q_n,\, k \neq k_0,k_1\}$ in non-decreasing
order of $|\sin 2 \pi (x+\frac {k p_n} {2 q_n})|$ (so $r=q_n-1$ or
$r=q_n-2$).  By (\ref {k-k''}), we have
$|\sin 2\pi(x+\frac{s_j p_n}{2q_n})|>C_1 \frac {j} {q_n}$.  Then
\begin{align}
\sum_{^{k=1}_{k \neq k_0,k_1}}^{q_n} \ln|\frac{\sin
  2\pi(x+\frac{k\al}{2})}{\sin
  2\pi(x+\frac{kp_n}{2q_n})}|&=
\sum_{1 \leq j\leq 4\frac {C_0} {C_1}} \ln|\frac{\sin
  2\pi(x+\frac{s_j \al}{2})}{\sin
  2\pi(x+\frac{s_j p_n}{2q_n})}|+\sum_{4\frac {C_0} {C_1}< j \leq r} \ln|\frac{\sin
  2\pi(x+\frac{s_j \al}{2})}{\sin
  2\pi(x+\frac{s_j p_n}{2q_n})}|\\
\nonumber
&\leq C \ln q_n+\sum_{4\frac {C_0} {C_1} < j \leq r}
\frac {C q_n \Delta_n} {j} \leq C \ln q_n,
\end{align}
which is (\ref {C ln q_n}).
\end{pf}

\comm{
$Let us estimate
\be
\sum_{^{k=1}_{k \neq k_0,k_1}}^{q_n} \ln|\frac{\sin
  2\pi(x+\frac{k\al}{2})}{\sin
  2\pi(x+\frac{kp_n}{2q_n})}|.
\ee

By (\ref{rat-1}) it is enough to estimate
$\displaystyle\sum_{^{k=1}_{k\not=k_0}}^{q_n}\ln|\frac{\sin
  2\pi(x+\frac{k\al}{2})}{\sin
  2\pi(x+\frac{kp_n}{2q_n})}|.$
We have
\be \la{16}
|\frac{\sin
  2\pi(x+\frac{k\al}{2})}{\sin
  2\pi(x+\frac{kp_n}{2q_n})}-1|=|\cos\frac{2\pi k\De_n}{2q_n}-1+
\frac{\cos
  2\pi(x+\frac{kp_n}{2q_n})}{\sin
  2\pi(x+\frac{kp_n}{2q_n})}\sin\frac{2\pi
  k\De_n}{2q_n}|<\frac{C\De_n}{|\sin
  2\pi(x+\frac{kp_n}{2q_n})|}
\ee
We can find $C_1<\infty$ so that if $|k-k_0|>C_1$ we have $\frac{C\De_n}{|\sin
  2\pi(x+\frac{kp_n}{2q_n})|}<\frac 14.$ Then,
arranging $k$ in the increasing order of $|\sin
  2\pi(x+\frac{kp_n}{2q_n})|,$ we obtain that 
\beq \la{17}
\big{|}\displaystyle\sum_{^{k=1}_{
k\not=k_0}}^{q_n}\ln|\frac{\sin
  2\pi(x+\frac{k\al}{2})}{\sin
  2\pi(x+\frac{kp_n}{2q_n})}|\big|&\le &\big{|}\displaystyle\sum_{^{|k-k_0|\le
  C_1}_{\;\;\;\;k\not=k_0}}\ln|\frac{\sin
  2\pi(x+\frac{k\al}{2})}{\sin
  2\pi(x+\frac{kp_n}{2q_n})}|\big|+\big{|}\displaystyle\sum_{|k-k_0|>C_1}\ln|\frac{\sin
  2\pi(x+\frac{k\al}{2})}{\sin
  2\pi(x+\frac{kp_n}{2q_n})}|\big|\nonumber\\&\le & C\ln q_n+
  Cq_n\De_n\sum_{i=1}^{\frac{q_n}2}\frac 1i\le C\ln q_n.
\eeq
}

\bl
Let $\ell \in \N$ be such that
$\ell<\frac {q_{r+1}} {10 q_n}$,
where $r\ge n$.
Given a sequence $|\ell_k|\le \ell-1,$
$\,k=1,\ldots,q_n,$ let $1\le k_0\le q_n$  be such that
\be
|\sin 2\pi(x+\frac{(k_0+ \ell_{k_0} q_r) \al}{2})|=
\min_{1\le k\le q_n}|\sin 2\pi(x+\frac{(k+\ell_k q_r)\al}{2})|.
\ee
Then
\be \la{shift1}
\left |\sum_{^{k=1}_{k\not=k_0}}^{q_n}\ln|\sin
  2\pi(x+\frac{(k+\ell_k q_r)\al}{2})|+(q_n-1)\ln 2 \right | < 
\ln q_n+C(\Delta_n+(\ell-1)\Delta_r) q_n\ln q_n.
\ee
%
\el

\begin{pf}

Notice that
$\|(k-k') \alpha\|_{\R/\Z} \geq \Delta_{n-1} \geq \frac {1} {2 q_n}$,
while $\|(\ell_k-\ell_{k'}) q_r \alpha\|_{\R/\Z}
\leq \frac {q_{r+1}} {5 q_n} \Delta_r<\frac {1} {5 q_n}$.  This implies
\be \label {k-k' 1}
\|(2x+(k+\ell_k q_r)\alpha)-(2x+(k'+\ell_{k'}q_r)\alpha)\|_{\R/\Z} \geq
\frac {1} {5 q_n}, \quad 1 \leq k,k' \leq q_n,k \neq k'.
\ee

\comm{
There are two cases to consider.  If $q_{n+1} \leq
q_n^{10}$ then we can use the simple bound (analogous to (\ref {k-k'}))
$\Delta_n$, which in this case is bounded by $q_n^{-C}$.
If $q_{n+1}>q_n^{10}$, then we have the
bound $|\ln \frac {\Delta_{n-1}} {2}| \leq C \ln q_n$ (since in this last
case we have $\ell \Delta_n \leq \frac {1} {40 q_n}$, so
$\|k-k_0 \alpha\|_{\R/\Z} \geq
\Delta_{n-1}$, but $\|(\ell_k-\ell_{k_0})q_r \alpha\|_{\R/\Z}
\leq \Delta_n/5$, which gives
$\|x+(k+\ell_k)q_r)\alpha\|_{\R/\Z}<\Delta_{n-1}/2$).
}

By (\ref {sin z-w}), we have
\be \label {sinka 1}
|\frac{\sin
  2\pi(x+\frac{(k+\ell_kq_r)\al}{2})}{\sin
  2\pi(x+\frac{kp_n}{2q_n})}-1| \leq \frac{C_0(\De_n+(\ell-1)\De_r)}{|\sin
  2\pi(x+\frac{kp_n}{2q_n})|}.
\ee
We now argue as in the previous lemma, using (\ref {k-k' 1}) and (\ref
{sinka 1}) instead of (\ref {k-k'}) and (\ref {sinka}).
\end{pf}




\subsection{Non-resonant case. Proof of Lemma \ref{non-res}}\la{nonres}
In the arguments that follow, we will actually consider a slightly larger
range of $k$, by assuming a weaker upper bound
$k \leq \max \{\frac {q_n} {20},50 q_{n+1}^{8/9}\}$.  The fact that the
estimates hold for this larger range will be useful later
(when dealing with the resonant case).

We start with the proof of the first part.  Let
$k=mq_n \pm (sq_{n-1}+r)=m q_n \pm k_0,\; s\ge 1,\;0\le r<q_{n-1}$, $k_0
\leq \frac {q_n} {2}$, be non-resonant.  Notice that $2 s q_{n-1}<q_n$. 
Assume without loss of generality that
$k=m q_n+k_0$, the other case being treated similarly.

Notice that if $m \geq 1$ then $k>\frac {q_n} {20}$ which implies that $k
\leq 50 q_{n+1}^{8/9}$, and we have
\be \label {m q_n+1}
m \leq \frac {50 q_{n+1}^{8/9}} {q_n}
\ee
(which is also obviously satisfied if $m=0$).

Set $I_1=[-[\frac {sq_{n-1}}2],sq_{n-1}-[\frac {sq_{n-1}}2]-1]$
and $I_2=[m q_n+k_0-[\frac {sq_{n-1}}2],mq_n+k_0+sq_{n-1}-
[\frac {sq_{n-1}}2]-1]$. Set
$\theta_j=\theta+j\alpha,j\in
I_1\cup I_2.$ The set $\{\theta_j\}_{j\in
I_1\cup I_2}$ consists of $2sq_{n-1}$ elements.
\bl \la{uni1}
For any $\eps>0$ and sufficiently
  large $n,$ set $\{\theta_j\}_{j\in
I_1\cup I_2}$  is $\frac{-2\ln\frac s{q_n}}{q_{n-1}}
+\eps$-uniform.
\el

\begin{pf}

We will first estimate the numerator in (\ref{prod}). We have

\begin{align} \la{num}
\sum_{^{j\in I_1\cup I_2}_{j\not=i}}&\ln |\cos 2\pi a-\cos
2\pi \theta_j|\\
\nonumber
&=\sum_{^{j\in I_1\cup I_2}_{j\not=i}}\ln
|\sin 2\pi\frac{a+\theta_j}2|+\displaystyle\sum_{^{j\in I_1\cup
    I_2}_{j\not=i}}\ln |\sin 2\pi\frac{a-\theta_j}2|+(2sq_{n-1}-1)\ln
2\\
\nonumber
&=\Sigma_+ + \Sigma_- + (2sq_{n-1}-1)\ln 2.
\end{align}

Both $\Sigma_+$ and $\Sigma_-$ consist of $2s$ terms
of the form of (\ref{irrat-1}) plus $2s$ terms of the form
\be
\ln \min_{j=1,\ldots,q_{n-1}}|\sin( 2\pi(x+\frac{j\al}{2}))|,
\ee
minus $\ln |\sin\frac{a\pm\theta_i}2|$. 
Therefore, by (\ref{irrat-1})

\be\la{num1}\displaystyle\sum_{^{j\in I_1\cup I_2}_{j\not=i}}\ln |\cos 2\pi a-\cos
2\pi\theta_j|\le -2sq_{n-1}\ln 2 + Cs\ln q_{n-1}.
\ee

To estimate the denominator of (\ref{prod}) we represent it again in
the form (\ref{num}) with $a= \theta_i.$ Assume that $i=j_0q_{n-1}+i_0\in
I_1,\; 0\le i_0<q_{n-1},$ the other
case being treated similarly.  Then
\be
\Sigma_-=\sum_{^{j\in I_1\cup I_2}_{j\not=i}}\ln |\sin \pi
(i-j)\al|.
\ee

On each interval $I\subset I_1$ of length $q_{n-1}$, minimum over
$t\in I$ of $|\sin \pi (t-i)\al|$ is achieved at $t-i$ of the form
$jq_{n-1}$ for some $j$.  This follows from the fact that if
$0<|z|<q_{n-1}$ and $2 |j| q_{n-1}<q_n$ then $\|(j q_{n-1}+z)
\alpha\|_{\R/\Z}>\|j q_{n-1} \alpha\|_{\R/\Z}$, since $\|z \alpha\|
\geq \Delta_{n-2}$ and $\|j q_{n-1} \alpha\|<\Delta_{n-2}/2$.  The possible
values of $j$ form an interval $[j^0_-,j^0_+]$ of size $s$
containing $j_0$.

Let now $T$ be an arbitrary interval of length $q_{n-1}$ contained in $I_2$.
Notice that $T$ is contained in $[i+m q_n+1,i+(m+1) q_n-1]$.
The minimum over $t\in T$ of $|\sin \pi (t-i)\al|$ is achieved at $t-i$ of
either the form $mq_n+jq_{n-1}$ or the form $(m+1)q_n-jq_{n-1}$ for some
$j\in \N$.\footnote{
Let $t \in T$ minimize $\|(t-i) \alpha\|_{\R/\Z}$, and let $j_u$,
$u \in \{0,1\}$ be such that $t_u=(m+u)q_n+(-1)^u j_u q_{n-1}+i \in T$.
If $t \neq t_0$ and $t
\neq t_1$ then $\|(t-t_u) \alpha\|_{\R/\Z} \geq \Delta_{n-2}$ as above.
Since the $(t_u-i) \alpha$ minus nearest integer are on opposite
sides of $0$, this implies that $\|(t_0-t_1) \alpha\|_{\R/\Z} \geq 2
\Delta_{n-2}$.  But one easily checks that $\|(t_0-t_1) \alpha\|_{\R/\Z}$ is
either equal to $\Delta_{n-2}$ (if $t_1>t_0$)
or to $\Delta_{n-1}+\Delta_{n-2}$ (if $t_1 \leq t_0$).}
For $u \in \{0,1\}$, let
$t_u \in T$ be (the unique number) of the form $t_u=i+(m+u) q_n+(-1)^u j_u q_{n-1}$ for some
$j_u \in \N$. Since $|t_u-t_{1-u}|<q_{n-1}$ it follows that 
\be \la {ju}
0\leq j_{1-u}+j_u-[\frac {q_n} {q_{n-1}}]\leq 1.
\ee
For all $j\in [1,[\frac {q_n} {q_{n-1}}]]$, we have the lower bound
\be
\|((-1)^u j q_{n-1}+(m+u)q_n)\alpha\|_{\R/\Z} \geq \Delta_{n-1}/2.
\ee
Indeed, by (\ref {m q_n+1}), if $m \geq 1$ then
$(m+u) \Delta_n \leq 100 \frac {q_{n+1}^{8/9}} {q_n} \Delta_n \leq
\Delta_{n-1}/2$, while $\|j q_{n-1} \alpha\|_{\R/\Z} \geq
\Delta_{n-1}$.  If $m=0$ then $(m+u)q_n+(-1)^u j q_{n-1} \in [1,q_n-1]$, and
we get the lower bound $\Delta_{n-1}$.  Those considerations also give the
upper bound
\be
(m+u)\Delta_n \leq \max \{\Delta_{n-1}/2,\Delta_n\}.
\ee
This gives the estimate, for all $j \in [1,[\frac {q_n} {q_{n-1}}]]$,
\be \la{ju1}
\|((-1)^u j q_{n-1}+(m+u)q_n)\alpha\|_{\R/\Z} \geq j \Delta_{n-1}/C.
\ee

Let $T$ now run through the set of disjoint segments
$T^p,$ each of length $q_{n-1}$, such that $I_2=\cup_{p=1}^sT^p.$ 
It is not difficult to see that there exists $u$ (possibly
both $u=0,1$) such that for all $p$ corresponding $j_u$ satisfy
$j_u \leq \frac {3} {4} [\frac {q_n} {q_{n-1}}]$.\footnote{
For $u=0,1$ the $j_u$ form an interval $[j_-^u,j_+^u]$
of length $s$ contained in $[1,[\frac {q_n} {q_{n-1}}]].$
If $j_+^u>\frac 34 [\frac {q_n} {q_{n-1}}],$ then,
since $s\leq [\frac {q_n} {2q_{n-1}}],$ we have
that $j_-^u>\frac 14 [\frac {q_n} {q_{n-1}}]+1.$  Then, by (\ref{ju}),
$j_+^{1-u}<\frac 34 [\frac {q_n} {q_{n-1}}].$}
We now fix $u\in\{0,1\}$ with this property.
Then $\|((m+u) q_n+(-1)^u j_u q_{n-1}) \alpha\|_{\R/\Z} \leq
\frac {3} {4}[\frac {q_n} {q_{n-1}}]\Delta_{n-1}+|m+u|\Delta_n\leq
(\frac {3} {4}[\frac {q_n} {q_{n-1}}]+1) \Delta_{n-1}.$
Then, by (\ref{ju}), $j_{1-u}\geq \frac 14  [\frac {q_n} {q_{n-1}}]$ and
by (\ref{ju1}), $\|((-1)^{1-u} j_{1-u} q_{n-1}+(m+1-u)q_n)\alpha\|_{\R/\Z}
 \geq \frac 1{4C}  [\frac {q_n} {q_{n-1}}] \Delta_{n-1}.$ Thus 
$\|(t_{1-u}-i) \alpha\|_{\R/\Z} \geq
\frac{\|(t_u-i) \alpha\|_{\R/\Z}}{7C}$.

Notice that the $j_u$
form an interval $[j_-,j_+]$ of length $s$, contained in $[1,[\frac
{3 q_n} {4 q_{n-1}}]]$.


Splitting again $\Sigma_-$ into $2s$ sums of length
$q_{n-1}$ and applying (\ref{irrat-1}) on each we obtain

\begin{align} \la{s-}
\Sigma_- > &-2sq_{n-1}\ln 2+
\sum_{^{j^0_- \leq j \leq j^0_+}_{j \neq j_0}}
\ln |\sin \pi (j-j_0)q_{n-1}\al|\\
\nonumber
&+\sum_{j=j_-}^{j_+} \ln
|\sin \pi ((-1)^u j q_{n-1}+(m+u)q_n)\al|-Cs-Cs\ln q_{n-1}.
\end{align}

Denote the sums in (\ref{s-}) by $\Sigma_1$ and $\Sigma_2.$
Since $j^0_- \leq j_0 \leq j^0_+$ and $|[j^0_-,j^0_+]|=s$
we have that
\be \la{s1}
\Sigma_1 > 2\sum_{j=1}^{[\frac{s} {2}]}\ln\sin|\pi jq_{n-1}\al|>
2\sum_{j=1}^{[\frac{s} {2}]}\ln 2j\Delta_{n-1}>s(\ln\frac s{q_n}-C).
\ee

For $j\in [j_-,j_+]$ we use (\ref{ju1}) to obtain

\be \la{s2}
\Sigma_2 > \sum_{j=1}^{s}\ln j\Delta_{n-1}-C s
>s\ln\frac s{q_n}-Cs.
\ee
Therefore,
\be\la{s-1}
\Sigma_- > -2sq_{n-1}\ln 2+2s(\ln\frac s{q_n}-
C\ln q_{n-1}).
\ee
$\Sigma_+$ is estimated in a similar way. 
Set
$J_1=[-[\frac{s+1}2],s-[\frac{s+1}2]-1]$ and
$J_2=[[s/2],s+[s/2]-1]$, which are two adjacent disjoint intervals of length
$s$.  Then $I_1\cup
I_2$ can be represented as a disjoint union of segments $B_j,\;j\in
J_1\cup J_2,$ each of length $q_{n-1}$. 
Applying (\ref{irrat-1}) on each $B_j$ we
obtain
\be \la{s+}
\Sigma_+ > -2sq_{n-1}\ln 2+
\sum_{j\in J_1\cup J_2}\ln
|\sin 2\pi\hat \theta_j|-Cs\ln q_{n-1}-\ln |\sin 2\pi (\theta+i\al)|
\ee
where
\be
|\sin 2\pi\hat \theta_j|=\min_{\ell \in B_j}|\sin 2\pi
(\theta + \frac{(i+\ell)\alpha}2)|.
\ee

Let $\tilde \theta_j=\hat \theta_j,\;j\in J_1$,
and $\tilde \theta_j=\hat \theta_j-\frac {m q_n \alpha} {2},\;j\in J_2$.
 Since
$\theta \notin \Theta$, for sufficiently large $n,$ we have that
$$
\min_{j\in J_1\cup J_2}|\sin 2\pi\tilde \theta_j|>\frac 1{9s^2q_{n-1}^2}.
$$ 
To estimate $|\sin 2\pi\hat \theta_j|,\;j\in J_2,$ we distinguish the two
cases:
\begin{enumerate}
\item If $q_{n+1}>(20 s^2q_{n-1}^2)^9,$ we write

\be\la{large}
|\sin 2\pi\hat \theta_j| \geq |\sin 2\pi\tilde \theta_j \cos \pi m\Delta_n|-
|\cos 2\pi\tilde \theta_j \sin \pi m\Delta_n|>
\frac 1{10 s^2q_{n-1}^2}-\frac 1{q_{n+1}^{1/9}}>\frac {1} {20 s^2q_{n-1}^2}.
\ee
\item If $q_{n+1}\le(20 s^2q_{n-1}^2)^9$ we use that  since
$\theta \notin \Theta,$ for large $n,$
\be \la{small}
\min_{j\in J_2}|\sin
  2\pi\hat \theta_j|>((2m+2)q_n)^{-2}>(4q_{n+1})^{-2}>(20 sq_{n-1})^{-36}
\ee
\end{enumerate}
In either case, 
$$
\ln \min_{j\in J_2}|\sin 2\pi\hat \theta_j|> -C\ln sq_{n-1}
$$

Let $J=J_1$ or $J=J_2$ and
assume that $\hat \theta_{j+1}=\hat \theta_j+\frac {q_{n-1}} {2} \alpha$ for
every $j,j+1 \in J$.  Applying again the Stirling formula we obtain
\be \label {J}
\sum_{j\in J}\ln
|\sin 2\pi\hat\theta_j|>-C\ln sq_{n-1} + \sum_{j=1}^{s}\ln \frac
{j\Delta_{n-1}} {C}>
s\ln\frac s{q_n}-C(\ln sq_{n-1}+s).
\ee

In the other case, decompose $J$ in maximal intervals $T_\kappa$
such that for $j,j+1 \in
T_\kappa$ we have $\hat \theta_{j+1}=\hat \theta_j+\frac {q_{n-1}} {2}
\alpha$.  Notice that the boundary points of an interval $T_\kappa$
are either boundary points of $J$ or satisfy
$\|2 \hat \theta_j\|_{\R/\Z}+\Delta_{n-1} \geq
\frac {\Delta_{n-2}} {2}$. 
Assuming $T_\kappa \neq J$, there exists $j \in
T_\kappa$ such that $\|2 \hat \theta_j\|_{\R/\Z}+\Delta_{n-1} \geq
\frac {\Delta_{n-2}} {2}$.
An estimate similar to (\ref {J}) gives
\be \label {JJ}
\sum_{j \in T_\kappa} \ln |\sin 2\pi\hat \theta_j|>-|T_\kappa|
\ln q_{n-1}-C(\ln sq_{n-1}+|T_\kappa|).
\ee
If $T_\kappa$ does not contain
a boundary point of $J$ (in particular $|T_k| \leq |J|-2=s-2$ and $[\frac
{q_n} {q_{n-1}}] \geq 2s \geq 6$), then
$T_\kappa$ does not contain any $j$ with $\|2 \hat
\theta_j\|_{\R/\Z}<\frac {\Delta_{n-2}} {10}<
\frac {\Delta_{n-2}} {2}-\Delta_{n-1}$ (otherwise
$|T_\kappa|-1 \geq \frac {\Delta_{n-2}} {\Delta_{n-1}}-2
\geq \frac {q_n} {2 q_{n-1}}-2 \geq
s-2$, which is impossible) and hence
\be \label {JJJJ}
\sum_{j \in T_\kappa} \ln |\sin 2\pi \hat\theta_j|>
-|T_\kappa| (\ln q_{n-1}+C).
\ee
Putting together all $T_{\kappa}$, using (\ref {JJ}) for the ones that
intersect the boundary of $J$ and (\ref {JJJJ}) for the others, we get
\be
\sum_{j\in J}\ln |\sin 2\pi\hat \theta_j|>s\ln\frac s{q_n}-C(\ln
sq_{n-1}+s)
\ee
in all cases.

Putting together $J=J_1$ and $J=J_2$ we have
\be
\sum_{j\in J_1\cup J_2}\ln
|\sin 2\pi \hat\theta_j|>
2s\ln\frac s{q_n}-C(\ln sq_{n-1}+s)
\ee

Combining it with (\ref{s+}) we obtain
\be \la{s+1}
\Sigma_+ > -2sq_{n-1}\ln 2+
2s(\ln\frac s{q_n}-C\ln q_{n-1})
\ee
Putting together (\ref{s+1}),(\ref{s-1}), and (\ref{num}) gives
\be
\displaystyle\sum_{^{j\in I_1\cup I_2}_{j\not=i}}\ln |\cos 2\pi \theta_i-\cos
2\pi \theta_j|>4s(\ln\frac {sq_{n-1}}{q_n}-C\ln q_{n-1}) - 2sq_{n-1}\ln 2
\ee
This together with (\ref{num1}) yields
$$
\max_{j\in I_1\cup I_2}\prod_{^{\ell\in I_1\cup I_2}_{\ell \not= j}}\frac{
|z-\cos 2\pi\theta_\ell|}
{|\cos 2\pi\theta_j-\cos 2\pi\theta_\ell|}<
e^{-4s\ln\frac {sq_{n-1}}{q_n}+Cs\ln q_{n-1}}
$$
as desired.
\end{pf}

By Lemmas \ref{nongen} and \ref{uni1} at least one
 of $\theta_j,\;j\in I_1\cup I_2,$ is not in $A_{2sq_{n-1}-1,L+
\frac{2\ln\frac{sq_{n-1}} {q_n}}{q_{n-1}}-\epsilon}$ where $\eps$
can be made arbitrarily small for large $n$.
By Lemma \ref{prop} and singularity
of $0$,\footnote {To get what we need here one can take in
Lemma \ref {prop}, besides $y=0$, also
$\epsilon=\frac {99} {100} L$ and $\delta=\frac {99} {400}$.}
we have that for all $j\in I_1,\;$  $\theta_j\in A_{2sq_{n-1}-1,L+
\frac{2\ln\frac {sq_{n-1}}{q_n}}{q_{n-1}}-\epsilon}$ (using that $(s+1)
q_{n-1}>q_n^{8/9}$ and the bound $\frac
{\ln q_n} {q_{n-1}}<L$).
Let $j_0\in I_2$ be such
 that $\theta_{j_0}\notin A_{2sq_{n-1}-1,L+
\frac{2\ln\frac {sq_{n-1}}{q_n}}{q_{n-1}}-\epsilon}.$
Set $I=[j_0-sq_{n-1}+1,j_0+sq_{n-1}-1]=[x_1,x_2].$ Then by
 (\ref{cra}),(\ref{pol0}),
\begin{align}
|G_I(k,x_i)|&<e^{(L+\eps_1)(2sq_{n-1}-2-|k-x_i|)-2sq_{n-1}(L+\frac{2\ln\frac
 {sq_{n-1}}{q_n}}{q_{n-1}}-\eps)}\\
\nonumber
&<e^{-(L+\eps_1)|k-x_i|-4sq_{n-1}\frac{\ln\frac
 {sq_{n-1}}{q_n}}{q_{n-1}}+(\eps_1+\eps)sq_{n-1}}.
\end{align} 
 Since \be\la{d/2}|k-x_i|\ge [\frac{sq_{n-1}}2]-1,\ee
 we
 obtain that 
\be\la{G1}|G_I(k,x_i)|<e^{-(L+9\frac{\ln\frac
 {sq_{n-1}}{q_n}}{q_{n-1}}-\epsilon)|k-x_i|}
\ee
 which in view of $(s+1)q_{n-1}>q_n^{8/9}$ gives the
 statement of the first part of Lemma \ref{non-res}.

We now assume $s=0$.  In this case $\al$ is ``Diophantine'' on the scale
$q_{n-1}$ however some caution is needed as it may not be so on the
scale $q_n$.  Let $k=m q_n \pm k_0,\; \max \{\frac {1} {20}
q_{n-1},q_n^{8/9}\}<k_0<q_{n-1}$.  We will assume that $m=q_n+k_0$, the
other case being analogous.

We distinguish three cases.

\begin{enumerate}

\item  If $\frac {1} {20} q_{n-1}< k_0 \le \frac {4} {5} q_{n-1}$, set
$I_1=[-[\frac{19}{40}q_{n-1}]+1,[\frac{19}{40}q_{n-1}]]$
and $I_2=[mq_n+[\frac{19}{40}q_{n-1}]+1,mq_n+2 [\frac {q_{n-1}} {2}]-
[\frac{19}{40}q_{n-1}]]$.

\item If $\frac {4} {5} q_{n-1}<k_0<q_{n-1}$ and $q_n \leq 2 q_{n-1}$,
define $I_1=[-[\frac {q_n} {4}]+1,[\frac {q_n} {4}]]$ and $I_2=[m q_n+[\frac
{q_n} {4}]+1,m q_n+2 [\frac {q_n} {2}]-[\frac {q_n} {4}]]$.

\item If $\frac {4} {5} q_{n-1}< k_0 < q_{n-1}$ and $q_n>2 q_{n-1}$, set
$I_1=[-[\frac{q_{n-1}}2]+1,q_{n-1}-[\frac{q_{n-1}}2]]$ and
$I_2=[mq_n+q_{n-1}-[\frac{q_{n-1}}2]+1,
mq_n+2q_{n-1}-[\frac{q_{n-1}} {2}]]$.

\end{enumerate}

Set $\theta_j=\theta+j\alpha,j\in
I_1\cup I_2.$ The set $\{\theta_j\}_{j\in
I_1\cup I_2}$ consists of $2[\frac{q_{n-1}}2]$ elements in the
first case, $2[\frac {q_n} {2}]$ elements in the second case,
and of $2q_{n-1}$ elements in the third case.

\bl \la{di}
For any $\eps>0$ and sufficiently
  large $n,$ the set $\{\theta_j\}_{j\in
I_1\cup I_2}$  is 
 $\eps$-uniform.
\el

\begin{pf}

Consider first the case $k_0 \le \frac {4} {5}
q_{n-1}$.  We will assume $q_{n-1}$ is even,
the other case needing obvious adjustments. As in the proof of Lemma
\ref{non-res} we will first 
estimate the numerator in (\ref{prod}).  We have

\begin{align} \la{num4}
\sum_{^{j\in I_1\cup I_2}_{j\not=i}}
\ln &|\cos 2\pi a-\cos
2\pi \theta_j|\\
\nonumber
&=\sum_{^{j\in I_1\cup I_2}_{j\not=i}}\ln
|\sin 2\pi\frac{a+\theta_j}2|+\sum_{^{j\in I_1\cup
    I_2}_{j\not=i}}\ln |\sin 2\pi\frac{a-\theta_j}2|+(q_{n-1}-1)\ln 2\\
\nonumber
&=\Sigma_+ +\Sigma_- + (q_{n-1}-1)\ln 2.
\end{align}

Both $\Sigma_+$ and $\Sigma_-$ are of the form
(\ref{shift1}) with $\ell_k \in \{0,m\}$\footnote {Recall that $m$ is chosen
so that $k=m q_n+k_0$, where $k \leq \max \{50 q_{n+1}^{8/9},q_n/20\}$. 
We have the bound $m \leq 50 q_{n+1}^{8/9}/q_n$, so (\ref {shift1}) really
applies.}
and $r=n$ plus a minimum term
minus $\ln |\sin 2\pi \frac{a\pm\theta_i}2|,$ so that the last two cancel
each other for the purpose of the upper bound.  Therefore, by (\ref{shift1})

\begin{align} \la{num5}
\sum_{^{j\in I_1\cup I_2}_{j\not=i}}\ln |\cos 2\pi a-\cos
2\pi\theta_j| &\le (1-q_{n-1}) \ln 2 +2\ln q_{n-1}+C(\De_{n-1}+
m\De_n)q_{n-1}\ln q_{n-1}\\
\nonumber
& \leq -q_{n-1} \ln 2+ C q_{n-1}^{8/9} \ln q_{n-1}.
\end{align}

To estimate the denominator of (\ref{prod}) we write it in
the form (\ref{num4}) with $a= \theta_i$. 
Then
\be
\Sigma_-=\sum_{^{j\in I_1\cup I_2}_{j\not=i}}\ln |\sin \pi
(i-j)\al|
\ee
is exactly of the form 
(\ref{shift1}).
 Therefore, by (\ref{shift1}),
\be \la{-0}
\Sigma_- >(1-q_{n-1})\ln 2-\ln q_{n-1}-C(\frac 1{q_{n}}+\frac 1{q_{n+1}^{1/9}})q_{n-1}\ln q_{n-1}>-q_{n-1}\ln2-Cq_{n-1}^{8/9}\ln q_{n-1}.
\ee
Similarly, for $\Sigma_+$ we have
\begin{align} \la{+0}
\Sigma_+ &>(1-q_{n-1})\ln2+\ln\displaystyle\min_{\ell \in I_1\cup I_2}
|\sin 2\pi (\theta + \frac{(i+\ell)\alpha}2)|-Cq_{n-1}^{8/9}\ln q_{n-1}\\
\nonumber
& >-q_{n-1}\ln2-Cq_{n-1}^{8/9}\ln q_{n-1}.
\end{align}
Here, we use the estimate 
\be\la{min+}
\ln \min_{\ell \in I_1\cup I_2} |\sin 2\pi (\theta +
\frac{(i+\ell)\alpha}2)|>-C\ln q_{n-1}
\ee
which is obtained by considering separately two cases $q_{n+1}>q_{n-1}^C$ and
$q_{n+1}<q_{n-1}^C,$ and arguing in the same way as in
(\ref{large}),(\ref{small}).
Combining (\ref{num5}),(\ref{num4}),(\ref{-0}) and  (\ref{+0}), we arrive at

\be \la{ep}
\max_{j\in I_1\cup I_2}\prod_{^{\ell\in I_1\cup I_2}_{\ell \not= j}}\frac{
|z-\cos 2\pi\theta_\ell)|}
{|\cos 2\pi\theta_j-\cos 2\pi\theta_\ell)|}<
e^{Cq_{n-1}^{8/9}\ln q_{n-1}}<e^{\eps q_{n-1}}
\ee
for any $\eps>0$ and sufficiently large $n,$ as stated.

For the other cases, $k>\frac {4} {5} q_{n-1}$, the proof is very similar.
If $q_n \leq 2 q_{n-1}$, the argument is the same (replacing $q_{n-1}$ by
$q_n$).  We will concentrate on the case $q_n>2 q_{n-1}$ where the changes
are slightly more substantial.
Arguing as above we obtain by (\ref{shift1})

\begin{align} \la{num6}
\sum_{^{j\in I_1\cup I_2}_{j\not=i}}\ln |\cos 2\pi a-\cos
2\pi\theta_j|&\le -2q_{n-1}\ln 2 +4\ln
q_{n-1}+C(\De_{n-1}+m\De_n)q_{n-1}\ln q_{n-1}\\
\nonumber &<-2q_{n-1}\ln 2+Cq_{n-1}^{8/9}\ln q_{n-1}
\end{align}
The denominator in (\ref{prod}) can be again split as $\Sigma_+ +
\Sigma_-+(2q_{n-1}-1)\ln2.$
Both $\Sigma_+$ and $\Sigma_-$ are,
up to a constant, the sums of two terms of the form
(\ref{shift1}) plus minimum terms
(two for $\Sigma_+$ and one for $\Sigma_-$).
For the minimum terms of $\Sigma_+$ the estimate
(\ref{min+}) holds so that we obtain
\be \la{+01}
\Sigma_+ >-2q_{n-1}\ln2-Cq_{n-1}^{8/9}\ln q_{n-1}.
\ee
For the minimum term of $\Sigma_-$,
that is, $\ln 
\min |\sin \pi (i-j)\al|$ (where the minimum is taken over all $j$ which
belong to the interval $I_1$ or $I_2$ that does not contain $i$)
we observe that it is achieved at $j_0$ such that
$\|(i-j_0) \alpha\|_{\R/\Z} \geq \Delta_{n-1}-m\Delta_n$ (since
the possible values of $|i-j|$ are
contained in $[m q_n+1,m q_n+2q_{n-1}-1]$ and $q_n>2 q_{n-1}$ by
hypothesis).
Thus, recalling that in the present situation we have $q_n^{8/9}<q_{n-1}$,
\be
\ln\min |\sin \pi (i-j)\al|> \ln(\De_{n-1}-m\De_n)>
\ln (\frac 1{2q_n}-\frac {50} {q_{n+1}^{1/9} q_n})>-C \ln q_n>
-C\ln q_{n-1}.
\ee
Therefore, by (\ref{shift1}),
\be \la{-01}
\Sigma_- >-2q_{n-1}\ln 2-Cq_{n-1}^{8/9}\ln q_{n-1}.
\ee
Combining  (\ref{num6}),(\ref{num4}),(\ref{-01}) and  (\ref{+01}), gives (\ref{ep}), as desired.
\end{pf}

By Lemmas \ref{nongen}
and \ref{di} at least one
 of $\theta_j,\;j\in I_1\cup I_2,$ is not in
$A_{2 [\frac {q_{n-1}} {2}]-1,L-\epsilon}$ if
$k_0 \le \frac {4} {5} q_{n-1},\;$
not in $A_{2 [\frac {q_n} {2}]-1,L-\epsilon}$ if $\frac {4} {5}
q_{n-1}<k_0<q_{n-1}$ and $q_n \leq 2 q_{n-1}$, and not in
$A_{2q_{n-1}-1,L-\epsilon}$ if $k_0 > \frac {4} {5}
q_{n-1}$ and $q_n>2q_{n-1}$, where $\eps$  can be made
arbitrarily small for large $n$.  By Lemma \ref{prop} and singularity
of $0$, we have that, in all three cases, for all
$j\in I_1$,  $\theta_j$ belongs to the corresponding $A_{\cdot,L-\epsilon}$.
Let $j_0\in I_2$ be such that $\theta_{j_0}\notin A_{L-\epsilon}.$

For $k_0 \le \frac {4} {5}
q_{n-1}$ set $I=[j_0-[\frac{q_{n-1}}2]+1,
j_0+[\frac{q_{n-1}}2]]=[x_1,x_2].$
We then have \be \la{d1}|k-x_i| > \frac{q_{n-1}}{40}.\ee  
Then by
 (\ref{cra}),(\ref{pol0}),
\be \la{G2}
|G_I(k,x_i)|<e^{(L+\eps_1)(q_{n-1}-2-|k-x_i|)-q_{n-1}(L-\eps)}<e^{-(L+\eps_1-40(\eps_1+\eps))|k-x_i|},
\ee
as desired. 

For $k_0>\frac {4} {5} q_{n-1}$ and $q_n \leq 2 q_{n-1}$, set $I=[j_0-[\frac
{q_n} {2}]+1,j_0+[\frac {q_n} {2}]]=[x_1,x_2]$.  Then
\be \label {d3}
|k-x_i|>\frac {q_n} {10},
\ee
since $k-x_1>\frac {4} {5}
q_{n-1}-\frac {q_n} {4} \geq \frac {3} {10} q_{n-1}$ and $x_2-k>\frac {3
q_n} {4}-q_{n-1}=\frac {3 q_{n-2}-q_{n-1}} {4}>\frac {q_{n-1}} {5}$ (using
that $\frac {4} {5} q_{n-1}<k_0 \leq \frac {q_n} {2}=\frac {q_{n-1}+q_{n-2}}
{2}$).  Thus for any $\eps>0$ and sufficiently large $n,$ by
(\ref{cra}),(\ref{pol0}), and estimating as in (\ref{G2})
\be \la{G3}
|G_I(k,x_i)|<e^{-(L-\eps)|k-x_i|}.
\ee

For $k_0>\frac {4} {5} q_{n-1}$ and $q_n>2q_{n-1}$ set
$I=[j_0-q_{n-1}+1,j_0+q_{n-1}-1]=[x_1,x_2]$.
Then
\be
\la{d2}|k-x_i| > \frac{3q_{n-1}}{10}.
\ee
This implies as before that (\ref {G3}) holds
for any $\epsilon>0$ and sufficiently large $n$.
This concludes the proof of Lemma \ref {non-res} in all cases.
\qed

The estimates in the proof of Lemma \ref {non-res} have the following
corollary which will be necessary later (when dealing with the resonant
case).

\bl \la{cor}
Fix $\eps>0$.  Assume $b_n<k \leq \max
\{\frac {q_n} {20},50 q_{n+1}^{8/9}\}$.
Let $d=\dist (k, \{\ell q_n\}_{\ell \ge 0})
>\frac 1{10}q_n$.  Let 
$\phi=\phi_E$ be a
generalised eigenfunction.  Assume that either
\begin{enumerate}
\item $q_n\ge q_{n-1}^{10/9}$, or
\item $q_n < q_{n-1}^{10/9}$ and $k<q_n^C$ for some $C<\infty$.
\end{enumerate}
Then, for sufficiently large $n$ ($n>n_0(\eps,c,E)$ in the first case,
$n>n_0(\epsilon,c,E,C)$ in the second case),
$$
|\phi(k)|<e^{-(L-\eps)\frac {d} {2}}.
$$
\el

\begin{pf}

Recall that the previous estimates in this subsection were obtained,
under the non-resonance hypothesis $\dist(k,\{\ell
q_n\}_{\ell \geq 0})>b_n$, for
$b_n<k \leq \max \{\frac {q_n} {20},50 q_{n+1}^{8/9}\}$.

If $q_n\ge q_{n-1}^{10/9}$, we have $s \geq [\frac {q_n^{1/10}} {10}]$, and
the statement follows immediately from
(\ref{poi}), (\ref{ge}) and (\ref{G1}), (\ref{d/2}).

In case $q_n<q_{n-1}^{10/9}$, (\ref{G1}), (\ref{d/2}), (\ref{G2}), (\ref{d1}),
 (\ref{G3}), (\ref {d3}), (\ref{d2}) only lead to $
|\phi(k)|<e^{-(L-\eps)cd}$ with certain $c<1/2$.  In order to prove the Lemma
as stated we will need an additional ``patching'' argument, which is very
similar to the one used in subsection \ref {pat}.


We will show that in this case
\be \la{cor1}
|\phi(k)|<e^{-(L-\eps)(d-\frac {q_{n-1}}{20})}
\ee
from which the statement of the lemma follows.
Assume $\ell q_n<k<(\ell+1)q_n$.
Using  (\ref{G1}), (\ref{d/2}), (\ref{G2}),(\ref{d1}) and 
 (\ref{G3}), (\ref {d3}), (\ref{d2}) we obtain that for every $y\in [\ell
  q_n,(\ell+1)q_n]$ with
$\dist(y,\{\ell q_n,(\ell+1)q_n\})>\frac
{1} {20}q_{n-1}$, there exists an interval $y\in I(y)=[x_1,x_2]\subset
[(\ell-1) q_n,(\ell+2)q_n]$

such that 
\be \la{I1}
\dist(y,\partial I(y))>\frac{q_{n-1}}{40},
\ee
\be \la {I12}
G_{I(y)}(y,x_i) < e^{-(L-\eps)|y-x_i|},\;i=1,2
\ee
(notice that under the condition $q_n<q_{n-1}^{10/9}$ we have
$b_n=\frac {q_{n-1}} {20}$).

We here denote the boundary of the interval $I(y)$, the set $\{x_1,x_2\},$ by
 $\partial I(y).$  For $z\in
\partial I(y)$ we let $z^\prime$ be the neighbor of $z,$ (i.e.,
 $|z-z^\prime|=1)$ not belonging to
 $ I(y).$ 

If $x_2+1<(\ell+1)q_n-\frac 1{20}q_{n-1},$ we expand
$\phi(x_2+1)$ in (\ref{poi}) iterating  (\ref{poi}) 
with $I= I(x_2+1),$ and if $x_1-1>\ell q_n+\frac
1{20}q_{n-1},$ we expand $\phi(x_1-1)$ in (\ref{poi}) iterating  (\ref{poi}) 
with $I=I(x_1-1).$
 We continue to expand each term of the
form $\phi(z)$  in the same fashion until we arrive to $z$ such that
either $z+1 \geq (\ell+1) q_n-\frac 1{20}q_{n-1}$, $z-1 \leq \ell q_n+\frac        
1{20}q_{n-1}$, or the number of $G_I$ terms in the product
becomes $[\frac{40d}{q_{n-1}}],$ whichever comes first.  We then obtain an
expression of the form
\be\la{block1}
\phi(k)=\displaystyle\sum_{s ; z_{i+1}\in\partial I(z_i^\prime)} 
G_{I(k)}(k,z_1) G_{I(z_1^\prime)}
(z_1^\prime,z_2)\cdots G_{I(z_s^\prime)}
(z_s^\prime,z_{s+1})\phi(z_{s+1}^\prime).
\ee
where in each term of the summation we have
$\ell q_n+\frac {1} {20} q_{n-1}+1<z_i<(\ell+1) q_n-\frac {1} {20}
q_{n-1}-1$, $i=1,\ldots,s,$ and 
either $z_{s+1} \notin [\ell q_n+\frac {1} {20} q_{n-1}+1,(\ell+1) q_n-\frac
{1} {20} q_{n-1}-1]$, $s+1 < [\frac{40d}{q_{n-1}}]$, or
$s+1= [\frac{40d}{q_{n-1}}]$.

By construction, for  each $z_i^\prime,\, i\le s,$ we have that
$I(z_i^\prime)$ is well-defined and satisfies (\ref{I1})
and (\ref {I12}).  We now 
consider the two cases,
$z_{s+1} \notin [\ell q_n+\frac {1} {20} q_{n-1}+1,(\ell+1) q_n-\frac
{1} {20} q_{n-1}-1]$, $s+1 < [\frac{40d}{q_{n-1}}]$, and
$s+1= [\frac{40d}{q_{n-1}}]$ separately. 
If $z_{s+1} \notin [\ell q_n+\frac {1} {20} q_{n-1}+1,(\ell+1) q_n-\frac
{1} {20} q_{n-1}-1]$, $s+1 < [\frac{40d}{q_{n-1}}]$,
we have, by (\ref{I12}) and (\ref{ge}),
\begin{align}
| G_{I(k)}(k,&z_1) G_{I(z_1^\prime)}
(z_1^\prime,z_2)\cdots G_{I(z_s^\prime)}
(z_s^\prime,z_{s+1})\phi(z_{s+1}^\prime)|\\
\nonumber
&\le e^{-(L-\eps)(|k-z_1|+\sum_{i=1}^{s}|z_i^\prime-z_{i+1}|)}(1+(\ell+2)
q_n)\\
\nonumber
&\le e^{-(L-\eps)(|k-z_{s+1}|-(s+1))}(1+(\ell+2) q_n)\le e^{-(L-\eps)(d-\frac{q_{n-1}}{20}-\frac{40d}{q_{n-1}})}
(1+q_{n-1}^C).
\end{align} 

If $s+1= [\frac{40d}{q_{n-1}}] ,$
using again  (\ref{ge}), (\ref {I12}), and also (\ref {I1}) we obtain
$$| G_{I(k)}(k,z_1) G_{I(z_1^\prime)}
(z_1^\prime,z_2)\cdots G_{I(z_s^\prime)}
(z_s^\prime,z_{s+1})\phi(z_{s+1}^\prime)|\le e^{-(L-\eps)\frac{q_{n-1}}{40}\frac{40d}{q_{n-1}}}(1+q_{n-1}^C).$$
 In either case,  
\be\la{est1}
 | G_{I(k)}(k,z_1) G_{I(z_1^\prime)}
(z_1^\prime,z_2)\cdots G_{I(z_s^\prime)}
(z_s^\prime,z_{s+1})\phi(z_{s+1}^\prime)|\le e^{-(L-2\eps)(d-\frac{q_{n-1}}{20})}
\ee
for $n$ sufficiently large. Finally, we observe that
the total
 number of terms in
 (\ref{block1})
is bounded above by $2^{[\frac{40d}{q_{n-1}}]}.$ Combining it with
(\ref{block1}),
(\ref{est1}) we obtain
$$
|\phi(k)|\le 2^{[\frac{40d}{q_{n-1}}]} e^{-(L-2\eps)(d-\frac{q_{n-1}}{20})}< e^{-(L-3\eps)(d-\frac{q_{n-1}}{20})}
$$      
for large $n$.
\end{pf}

\subsection{Resonant case. Proof of Lemma \ref{res}}\la{ress}

Notice that, under the condition that $k$ is resonant, we have $k \geq \frac
{q_n} {2}$, which implies $q_{n+1}^{8/9} \geq \frac {q_n} {2}$.  This
is an implicit hypothesis in the next lemma.

\bl\la{zero}For any $\eps>0,$ for sufficiently large $n,$
and any $b\in [-\frac{13}8q_n,-\frac 38q_n]\cap\bbz,$ we have
$\theta+(b+q_n-1)\al\in A_{2q_n-1,\frac{23L}{32}+\eps}.$
\el
\begin{pf}

Let $b_1=b-1,\;b_2=b+2q_n-1.$ 

Applying Lemma \ref{cor} 
we obtain that for $i=1,2,$
\be
|\phi_E(b_i)|<
\left\{\;
\begin{array}{cc}
& e^{-(L-\eps)(b/2+q_n)},\;\; -\frac{13}8q_n\le b\le -
\frac{3}2q_n, \\[0.15 in]
& e^{-(L-\eps)|\frac{b+q_n}2|}
,\;\;-\frac{3}2q_n \le b\le -\frac{q_n}2,\;|b+q_n|>\frac{q_n}4,  \\[0.15 in]
& e^{(L-\eps)\frac b2},\;\;-\frac{q_n}2\le b\le
-\frac{3}8q_n.
\end{array} \right.
\ee

Using (\ref{poi}) with $I=[b,b+2q_n-2]$ we get
\be \la{G2'}
\max (|G_I(0,b)|,|G_I(0,b+2q_n-2)|)>\left\{\;
\begin{array}{cc}
& e^{(L-\eps)(b/2+q_n)},\;\; -\frac{13}8q_n\le b\le -\frac{3}2q_n, \\[0.15 in]
& e^{(L-\eps)|\frac{b+q_n}2|}
,\;\;-\frac{3}2q_n \le b\le -\frac{q_n}2,\;|b+q_n|> \frac{q_n}4, \\[0.15 in]
& e^{-(L-\eps)\frac b2},\;\;-\frac{q_n}2\le b\le -\frac{3}8q_n,\\[0.15 in]
& e^{-\eps q_n},\;\;|b+q_n|<\frac{q_n}4.

\end{array} \right.
\ee
By (\ref{cra}),(\ref{pol0}),

\begin{align} \la{G4}
&|Q_{2q_n-1} (\cos 2\pi(\theta+(b+q_n-1)\al))|\\
\nonumber
&=|P_{2q_n-1}(\theta+b\al)|<\min \{|G_I(0,b)|^{-1}
e^{(L+\eps_1)(b+2q_n-2)},|G_I(0,b+2q_n-2)|^{-1}e^{-
(L+\eps_1)b}\}.
\end{align}

Therefore, using (\ref{G2'}),(\ref{G4}) we obtain that
$\theta+(b+q_n-1)\al$ belongs to

\begin{itemize}\item $A_{2q_n-1,\frac{23L}{32}+\eps},$ if $-\frac{13}8q_n\le b\le -\frac{3}2q_n$ or $-\frac{q_n}2\le b\le
-\frac{3}8q_n, $ 

\item
$A_{2q_n-1,\frac{5L}{8}+\eps},$ if $-\frac{3}2q_n \le b\le -\frac{q_n}2.$


\end{itemize}
for any $\eps>0$ and sufficiently large $n.$
\end{pf}

 Fix $ 1\le\ell\le q_{n+1}^{8/9}/q_n$.  Set
$I_1=[-[\frac{5}8q_n],[\frac 58q_n]-1]$ and 
$I_2=[(\ell -1) q_n+
[\frac{5}8q_n],(\ell+1) q_n-[\frac 58q_n]-1].$ Set
  $\theta_j=\theta+j\al,j\in I_1\cup I_2.$ 
\bl \la{final}
Assume $L>\frac {16} {9} \beta.$ There exists an $\eps>0$ such that for
 sufficiently large $n,$ set $\{\theta_j,\,j\in I_1\cup I_2\}$ is 
$(\frac {9L}{32}-\eps)$-uniform.
\el 

We will now finish the proof of Lemma \ref{res} and prove Lemma
\ref{final} at the end of the subsection.

Let $k$ be resonant. Assume without loss of generality that $k=\ell
q_n+r$, $0\le r\le \max \{q_n^{8/9},\frac{q_{n-1}}{20}\}$,
$1\le\ell\le q_{n+1}^{8/9}/q_n$.

By Lemmas \ref{nongen},\ref{zero},\ref{final} there is $j_0\in I_2 $
such that $\theta+j_0\al\notin  
A_{2q_n-1,\frac{23L}{32}+\eps}.$ Set $I=[j_0-q_{n}+1,j_0+q_{n}-1]=[x_1,x_2].$
Then $$|G_I(k,x_i)|< 
e^{(L+\eps_1)(2q_{n}-2-|k-x_i|)-2q_{n}
(\frac{23L}{32}+\eps)}<e^{q_n(\frac{9L}{16}+\eps)
-(L+\eps)|k-x_i|}.$$ 
 Since, by a simple computation, $|k-x_i|>(5/8-\eps-\frac{1}{20})q_{n},$ we
 obtain that 
\be \la{G3'}
|G_I(k,x_i)|<e^{(-\frac L{46}+\eps)|k-x_i|}
\ee
 which gives the
 statement of Lemma \ref{non-res}. 
\qed

\smallskip
\noindent {\it Proof of Lemma \ref{final}.} 
As in the proof of Lemma \ref{non-res} we will first 
estimate the numerator in (\ref{prod}). We have

\begin{align} \la{num2}
\sum_{^{j\in I_1\cup I_2}_{j\not=i}} &\ln |\cos 2\pi a-\cos
2\pi \theta_j|\\
\nonumber
&=\sum_{^{j\in I_1\cup I_2}_{j\not=i}}\ln
|\sin 2\pi\frac{a+\theta_j}2|
+\displaystyle\sum_{^{j\in I_1\cup
    I_2}_{j\not=i}}\ln |\sin 2\pi\frac{a-\theta_j}2|+(2q_{n}-1)\ln 2
\\
\nonumber
&=\Sigma_+ + \Sigma_- + (2q_{n}-1)\ln 2.
\end{align}

Both $\Sigma_+$ and $\Sigma_-$ consist of $2$ terms
of the form of (\ref{shift1}) with $r=n,$ plus two terms of the form $\ln
\min_{k=1,\ldots,q_{n}}|\sin 2\pi(x+\frac{(k+\ell_kq_n)\al}{2})|,$
where $\ell_k\in\{0,\pm(\ell-1),\pm \ell\},\;
k=1,\ldots,q_{n},$ minus $\ln
|\sin2\pi\frac{a\pm\theta_i}2|.$ 
%

Therefore, by (\ref{shift1})

\be\la{num3}
\sum_{^{j\in I_1\cup I_2}_{j\not=i}}\ln |\cos 2\pi a-\cos
2\pi\theta_j|\le (2-2q_{n})\ln 2 +4\ln q_n+C\ell \De_nq_n\ln q_n.
\ee

To estimate the denominator of (\ref{prod}) we write it in
the form (\ref{num2}) with $a= \theta_i.$ 
Then $\Sigma_-
=\displaystyle\sum_{^{j\in I_1\cup I_2}_{j\not=i}}\ln |\sin \pi
(i-j)\al|$ can be split into two sums of the form 
(\ref{shift1}) plus the
minimum term. The corresponding minimum term is
achieved at $|i-j_0|$ of the form $q_n$ or $\ell q_n$.
Therefore, for any
$\eps_1>0$ and sufficiently large $n$

\begin{align} \la{s-res}
\Sigma_-&>-2q_n\ln2+ \ln |\sin \pi
q_n\al|-C\max
(\ln q_n,\ell \De_nq_n\ln q_n)\\
\nonumber
&>-2q_n\ln2-\ln q_{n+1}-C\max
(\ln q_n,\ell \De_nq_n\ln q_n).
\end{align}

Since
$$
\sin 2\pi(\theta+\frac{(k+i+\ell_kq_n)\al}2)=\sin
2\pi(\theta+\frac{(k+i)\al}2)\cos \pi \ell_k\Delta_n
\pm \cos 2\pi(\theta+\frac{(k+i)\al}2)\sin\pi\ell_k\Delta_n
$$
(the $\pm$ depending on the sign of $q_n \alpha-p_n$)
we have, by (\ref{Theta}) that if $q_{n+1} \geq q_n^{10}$ then
\be
\min_{k,i\in [-q_n,q_n-1],\,\ell_k\in\{0,\pm(\ell-1),\pm\ell\}} |\sin
2\pi(\theta+\frac{(k+i+\ell_kq_n)\al}2)|>\frac {1} {10} q_n^{-2},
\ee
and if $q_{n+1}<q_n^{10}$ then we have the obvious
\be
\min_{k,i\in [-q_n,q_n-1],\,\ell_k\in\{0,\pm(\ell-1),\pm\ell\}} |\sin
2\pi(\theta+\frac{(k+i+\ell_kq_n)\al}2)|>\frac 15 q_{n+1}^{-2}>\frac {1} {5}
q_n^{-20}.
\ee

As before, $\Sigma_+$ can be split into two sums of the form
(\ref{shift1}) plus two minimum terms minus $\ln|\sin 2\pi(\theta+i\al)|.$
Therefore,
\be\la{s+res}
\Sigma_+ >-2q_n\ln2-C\max
(\ln q_n,\ell \De_nq_n\ln q_n)
\ee
Combining (\ref {num2}),(\ref{num3}),(\ref{s-res}) and (\ref{s+res})
we obtain

\be\la{fin}
\max_{j\in I_1\cup I_2}\prod_{^{\ell\in I_1\cup I_2}_{\ell \not= j}}\frac{
|z-\cos 2\pi\theta_\ell)|}
{|\cos 2\pi\theta_j-\cos 2\pi\theta_\ell)|}<q_n^Ce^{(\beta+\eps_1)q_n}
\ee

For $\beta <\frac 9{16}L$ this gives the desired bound.\qed

\end{document}